\documentclass[11pt,reqno]{amsart}

\usepackage{amscd}
\usepackage{amsfonts}
\usepackage{amsmath}
\usepackage{amssymb}
\usepackage{amsthm}
\usepackage{fancyhdr}
\usepackage{latexsym}
\usepackage[colorlinks=true, pdfstartview=FitV, linkcolor=blue,
            citecolor=blue, urlcolor=blue]{hyperref}

\usepackage{cancel}

\input epsf
\input texdraw
\input txdtools.tex
\input xy
\xyoption{all}

\usepackage{color}

\definecolor{Red}{rgb}{1.00, 0.00, 0.00}
\definecolor{DarkGreen}{rgb}{0.00, 1.00, 0.00}
\definecolor{Blue}{rgb}{0.00, 0.00, 1.00}
\definecolor{Cyan}{rgb}{0.00, 1.00, 1.00}
\definecolor{Magenta}{rgb}{1.00, 0.00, 1.00}
\definecolor{DeepSkyBlue}{rgb}{0.00, 0.75, 1.00}
\definecolor{DarkGreen}{rgb}{0.00, 0.39, 0.00}
\definecolor{SpringGreen}{rgb}{0.00, 1.00, 0.50}
\definecolor{DarkOrange}{rgb}{1.00, 0.55, 0.00}
\definecolor{OrangeRed}{rgb}{1.00, 0.27, 0.00}
\definecolor{DeepPink}{rgb}{1.00, 0.08, 0.57}
\definecolor{DarkViolet}{rgb}{0.58, 0.00, 0.82}
\definecolor{SaddleBrown}{rgb}{0.54, 0.27, 0.07}
\definecolor{Black}{rgb}{0.00, 0.00, 0.00}
\definecolor{dark-magenta}{rgb}{.5,0,.5}
\definecolor{myblack}{rgb}{0,0,0}
\definecolor{darkgray}{gray}{0.5}
\definecolor{lightgray}{gray}{0.75}

\usepackage[pdftex]{graphicx}
\usepackage{epstopdf}

\newcommand{\medno}{\medskip\noindent}

\newcommand{\rr}{\mathbb{R}}
\newcommand{\p}{\partial}
\newcommand{\zz}{\mathbb{Z}}

\newcommand{\ci}{\mathbb{T}}

\newcommand{\tor}{\mathbb{T}}
\newcommand{\ee}{\varepsilon}

\theoremstyle{plain}  
\newtheorem{theorem}{Theorem}
\newtheorem{proposition}{Proposition}
\newtheorem{lemma}{Lemma}
\newtheorem{corollary}{Corollary}

\theoremstyle{definition}

\usepackage{tikz}

\usetikzlibrary{patterns}

\begin{document}

\title{
Construction of 2-peakon Solutions \\ and Ill-Posedness
for  the Novikov equation
}
  
\author{ A. Alexandrou Himonas,  Curtis Holliman  \& Carlos Kenig}



\begin{abstract}
For the Novikov equation, on both the line and  the circle,
we construct  a 2-peakon solution 
with an asymmetric antipeakon-peakon initial profile
whose  $H^s$-norm for $s<3/2$ is  arbitrarily small. 
Immediately after the initial time,  both the antipeakon and peakon  move  in the positive direction, 
and
a collision occurs in arbitrarily small time. Moreover,
at the collision time the $H^s$-norm of the solution becomes arbitrarily large  when $5/4<s<3/2$, thus resulting in norm inflation and  ill-posedness.
However, when $s<5/4$, 
the solution at the collision time coincides with a second solitary antipeakon solution. 
This scenario thus results in 
nonuniqueness and ill-posedness. Finally, when $s=5/4$ ill-posedness follows
either from a failure of convergence or a failure of uniqueness. Considering that
the Novikov equation is well-posed for $s>3/2$, these results put together
establish $3/2$ as the critical index of well-posedness for this equation. The case $s=3/2$ remains an open question.
\end{abstract}

\date{August 18, 2017;\, Corresponding author: himonas.1@nd.edu}

\maketitle

\markboth{}
{}



%
%
%
%
%
%
%
\section{Introduction and Results}
\label{sec:1}
 \setcounter{equation}{0}
We consider the Cauchy problem
for the Novikov equation (NE) on the line and the circle
\begin{align}
 \label{NE}
&u_t  + u^2u_x
 +
 \p_x D^{-2}
\Big[
 u^3+ \frac{3}{2} uu_x^2
\Big]
+
D^{-2}
\Big[
  \frac{1}{2}u_x^3 
\Big]=0, 
\\&
u(x,0) = u_0(x), \ \ x\in \rr \text{ or } \ci, \ \ t\in\rr,
 \label{NE-data}
\end{align}
where $D^{-2}$ is the Bessel potential $D^{-2}=(1- \p_x^2)^{-1}$, 
and construct specific 2-peakon solutions $u(t)$ that 
collide at a finite time $T$ in such a way as to give rise to the phenomenon
of norm-inflation.  In particular, the norm-inflation generated by these 2-peakon collisions
occur 
in Sobolev spaces $H^s$ with exponents between 5/4 and 3/2.  As such, we will
refer to 3/2 as the {\it critical} exponent for well-posedness, as well-posedness
has been proven for exponents greater than 3/2 (see \cite{hh2}).  
For exponents $s$ less than $5/4$, the collision of the 2-peakons in fact converges to a single antipeakon $u(T)$, which can be thought of as a superposition of  both peakons.
This scenario allows us to demonstrate non-uniqueness.
Taken together, these results prove that NE is ill-posed in $H^s$ for $s<3/2$.

\vspace{.1in}
We recall that  NE is  well-posed  in the sense of 
Hadamard (see \cite{h})  in Sobolev spaces $H^s$ with exponents $s > 3/2$  (see \cite{hh2}). More precisely,
if $u_0$ belongs to the Sobolev space $H^s$ on the circle or the line, then there exists $T_s=T_s(||u_0||_{H^s})>0$ and a unique solution $u \in C([0,T_s];H^s)$ of the Cauchy problem for the Novikov equation \eqref{NE}--\eqref{NE-data} satisfying the following  estimate
\begin{equation}
\label{gCH-size-est}
\|u(t)\|_{H^s} \leq 2\|u_0\|_{H^s}, \ \ 
\text{for} \ \ 
0 \leq t \leq T_s,
\quad
\text{with}
\quad
 T_s
=
   \frac{1}{4c_s\|u_0\|_{H^s}^2},
\end{equation}
where  $c_s>0$ is a constant depending on $s$.
 Furthermore, the data-to-solution map $u(0)\mapsto u(t)$ is 
 continuous but not uniformly continuous.

 %
 %
 %
 %
 \vspace{.1in}
 The Novikov equation   is an integrable equation and  its local form,
\begin{equation}
\label{ne}
(1-\p_x^2)u_t=u^2u_{xxx}+3uu_xu_{xx}-4u^2u_x,
\end{equation}
was derived  by Vladimir Novikov \cite{nov}
in his attempt to classify all integrable Camassa-Holm--type equations
with quadratic and cubic nonlinearities  of the form 
$
(1-\p_x^2)u_t=P(u,u_x,u_{xx},...),
$
where $P$ is a polynomial of $u$ and its derivatives.  
The Lax pair for NE was derived by Hone and Wang in \cite{hw2008}
and is given by the equations
\begin{subequations}\label{laxpair-ne}
\begin{equation}
\left(
\arraycolsep=1pt\def\arraystretch{1}
\begin{array}{c}
\psi_1
\\
\psi_2
\\
\psi_3
\end{array}
\right)_{\!x}
=
U\left(m, \lambda\right)
\left(
\arraycolsep=1pt\def\arraystretch{1}
\begin{array}{c}
\psi_1
\\
\psi_2
\\
\psi_3
\end{array}
\right),
\quad
\left(
\arraycolsep=1pt\def\arraystretch{1}
\begin{array}{c}
\psi_1
\\
\psi_2
\\
\psi_3
\end{array}
\right)_{\!t}
=
V\left(m, u, \lambda\right)
\left(
\arraycolsep=1pt\def\arraystretch{1}
\begin{array}{c}
\psi_1
\\
\psi_2
\\
\psi_3
\end{array}
\right),
\end{equation}
where  $m=u-u_{xx}$ and  the matrices $U$ and $V$ are defined by
\begin{equation}
U\left(m, \lambda\right)
=\!
\left(
\arraycolsep=3pt\def\arraystretch{1}
\begin{array}{lcr}
0 & \lambda m & 1
\\
0 &0 & \lambda m  
\\
1 &0 &0
\end{array}
\right)\!,
\,
V\left(m, u, \lambda\right)
=\!
\left(
\!
\arraycolsep=4pt\def\arraystretch{1.2}
\begin{array}{lcr}
\frac{1}{3\lambda^2}-uu_x
& \frac{u_x}{\lambda}-\lambda m u^2 
& u_x^2
\\
\frac{u}{\lambda}
&-\frac{2}{3\lambda^2}
&-\frac{u_x}{\lambda}-\lambda m u^2
\\
-u^2 & \frac{u}{\lambda} &\frac{1}{3\lambda^2}+uu_x
\end{array}
\!
\right)\!.
\end{equation}
\end{subequations}
The Novikov equation possesses peakon traveling wave solutions 
   \cite{hm},\cite{hls}, \cite{gh}, which on the real line are  given by the formula
\begin{equation}
\label{NE-peak}
u(x,t)=\pm \sqrt{c}\, e^{-|x-ct|},
\end{equation}
where $c>0$ is the wave speed.  On the circle, the peakon solutions
are given by the formula
\begin{equation}
\label{NE-peak-per}
u(x,t)=\frac{\sqrt{c}}{\cosh(\pi)}\cosh([x-ct]_p-\pi),
\quad
\text{where}
\quad
[x-ct]_p \doteq x-ct-2\pi\Big[\frac{x-ct}{2\pi}\Big].
\end{equation}
In fact, the Novikov equation possess multi-peakon traveling 
wave solutions on both the line and the circle  \cite{hm},\cite{hls}, \cite{gh}. 
More precisely, on the line the  $n$-peakon,
\begin{equation}
\label{mpeak}
u(x,t) =\sum_{j=1}^n p_j(t) e^{-|x-q_j(t)|},
\end{equation}
is a solution  to NE
if and only if the positions $(q_1,\cdots, q_n)$ and  the momenta
$(p_1,\cdots, p_n)$ satisfy the following system of  $2n$ differential equations 
\begin{equation}
\label{odes}
\begin{cases}
\,\,\displaystyle  \frac{dq_j}{dt}&=u^2(q_j),
\\[9pt]
\,\,\displaystyle  \frac{dp_j}{dt}&=-u(q_j)u_x(q_j)p_j .
\end{cases}
\end{equation}
The description of the periodic $n$-peakon is similar.
Furthermore,  NE solutions conserve the  $H^1$-norm, that is 
\begin{align}
\label{h1-conserv}
\int_{\rr  \text{ or } \ci}
 \Big[u^2(t) + u_x^2(t)\Big] \, dx
 =
 \int_{\rr  \text{ or } \ci}
 \Big[u^2(0) + u_x^2(0)\Big] \, dx.
\end{align}

Next, we state our first result that gives the basic properties of the 2-peakon solutions, which are constructed here and are needed 
for proving  the ill-posedness of NE below $3/2$.

\begin{theorem}
\label{NE-2peak-slns}
For any $\ee > 0$ there exists a $T > 0$ for which the NE Cauchy problem on the line and the circle
 \eqref{NE}--\eqref{NE-data} 
has a 2-peakon solution $u \in C([0,T];H^s)$ such that its lifespan  and 
its initial size  satisfy the estimates
\begin{align}
\label{norm-inflation-1}
&\text{Lifespan} =T < \ee, \\
\label{norm-inflation-2}
&\|u_0\|_{H^s} < \ee,
\end{align}
while as  $t$ approaches the lifespan $T$  the $H^s$ norm
of the solution $u(t)$ satisfies the estimates
\begin{equation}
\label{limits}
\lim_{t\to T} \|u(t)\|_{H^s}
=
\begin{cases}
\infty\quad
\text{(norm inflation)},     \hskip0.5in    5/4 < s < 3/2, \\
\text{may not exist},  \hskip1.0in  s = 5/4, \\
C_s,  \quad \text{for some $C_s>0$}, \hskip0.45in  s < 5/4.
\end{cases}
\end{equation}
Moreover, when $s<5/4$ then $u(t)$ converges 
to an antipeakon
 $u(T)= -\sqrt{c_T} \, e^{-|x - q_T|}$, for some $c_T>0$ and 
 $q_T>0$, with   $\|u(T)\|_{H^s}=C_s$.
\end{theorem}
This theorem is a very interesting result in its own right.  
Unlike the  Camassa-Holm (CH) equation
(see  \cite{cam},  \cite{ff},  \cite{l1}, \cite{mn})
\begin{align}
\label{ch-eq}
(1-\p_x^2)u_t=uu_{xxx}+2u_xu_{xx}-3uu_x
\end{align}
and the Degasperis-Procesi (DP) equation 
(see \cite{dp}, \cite{hs}, \cite{ls},
\cite{l2}, \cite{dhh})
\begin{align}
\label{dp-eq}
(1-\p_x^2)u_t=uu_{xxx}+3u_xu_{xx}-4uu_x,
\end{align}
for which we can construct  special symmetric 
2-peakon solutions, called  peakon-antipeakons, 
of the form
\begin{align}
\label{pea-antipeak}
u(x,t) 
=
p(t) e^{-|x+q(t)|} - p(t) e^{-|x-q(t)|},
\end{align}
this is impossible for NE.
Peakon-antipeakon solutions, which are  convenient to work with, are possible for CH and DP because these equations
contain a symmetry that allows us to reduce 
the corresponding to \eqref{odes} ODE system for the positions and the
momenta  via 
$p = p_1 = -p_2$ and $q = q_1 = -q_2$.
This symmetry causes the peak and antipeak to move against 
each other and collide
in finite time (see \cite{hhg}, \cite{hgh}, \cite{by}).  
Such a construction is not possible 
for NE because  by equations \eqref{odes}
we have $\frac{dq_j}{dt}\ge 0$ for all positions $q_j$.
Thus, we see that for NE all the peaks and antipeaks move in the same direction. 
Therefore collision can occur only if
the peakon that  follows moves faster than the one ahead of it, and eventually overtakes it. For this scenario to happen we must 
break symmetry and solve 
the full system of the four highly nonlinear 
differential equation defined  by system \eqref{odes}
for $n=2$ with appropriate initial data.
This procedure involves several novel ideas 
which are described in the Sections \ref{sec:2} and  \ref{sec:7}.
The results are summarized
 in Theorems  \ref{cor:NE-2-peakon} and  
 \ref{per-cor:NE-2-peakon}.

\vskip0.1in
Next, using  Theorem \ref{NE-2peak-slns}
we obtain the following ill-posedness result for NE.

\begin{theorem}
\label{NE-ILP}
The Cauchy problem  for the Novikov equation
on the line and the circle
 \eqref{NE}--\eqref{NE-data} 
is ill-posed in Sobolev spaces $H^s$ for $s<3/2$.
More precisely, if $5/4<s<3/2$ then the data-to-solution map
is not continuous while  if $s<5/4$ then solution is not unique. 
When $s=5/4$ then either continuity or  uniqueness fails.
\end{theorem}
As we have mentioned before,  this theorem
combined with the well-posedness result 
for NE in  $H^s$, $s>3/2$, 
proved in  \cite{hh2}, completes 
the well-posedness picture of NE in Sobolev spaces, 
except for $s=3/2$, which remains an open question.
It is worth comparing the ill-posedness of NE, 
which has cubic nonlinearities, 
with those  of CH and  DP, 
the two integrable equations of the same type 
but with quadratic nonlinearities, 
which are both well-posed in $H^s$, $s>3/2$.
Defining  the {\it ``inflation index"} to be the 
Sobolev exponent  $s_i$ such that there is norm inflation
(which implies discontinuity of the data-to-solution map)
for all $s_i<s<3/2$,
we have the following observations.
For  CH  the inflation index $s_i=1$ and coincides with 
the index of the $H^1$-norm, which is the most important 
 conserved quantity  of CH. For $s<1$ 
 the peakon-antipeakon traveling wave solution 
\eqref{pea-antipeak}  for CH
converges  in $H^s$ to $u(T)=0$ as $t$ approaches the collision time $T$, giving rise to another solution 
(namely the trivial solution) demonstrating 
 ill-posedness due to failure of uniqueness.
However, for  DP the inflation index $s_i=1/2$.
 When  $s<1/2$, then the corresponding peakon-antipeakon traveling wave solution for DP converges  in $H^s$ to a function,
 which gives rise to another kind of DP solution,
 called ``shock peakon" that results
 to failure of uniqueness (see \cite{hhg}).
 From our results above we see that the 
 the inflation index for NE is  $5/4$, which is a very 
 interesting number and which 
 follows from the limiting behavior of the momenta $p_1(t)$
 and $p_2(t)$ as $t$ approaches the collision time
 (see  Theorem \ref{cor:NE-2-peakon}).
 For $s<5/4$ it is shown that the  2-peakon solution
\eqref{2-peak} constructed in Section  \ref{sec:2}
converges  in $H^s$ to an antipeakon, 
 which gives rise to an antipeakon traveling 
 wave solution demonstrating failure of uniqueness
 (see Proposition  \ref{nu-hs}).

\vskip0.1in
 Finally, we mention that  the  method used here for  proving ill-posedness for  NE is similar to that used by many authors for other nonlinear evolution equations.
For example, Bourgain and Pavlovic in  \cite{bp} proved  ill-posedness 
for the 3D Navier-Stokes equations in  Besov spaces in the sense of norm inflation.
Similar  methods  for establishing ill-posedness for dispersive equations 
have been used  by  Kenig, Ponce and Vega \cite{kpv} and Christ, Colliander and  Tao  \cite{cct}. 
The ill-posedness for the generalized KdV and nonlinear Schr\"{o}dinger equations in Sobolev spaces
has been tackled in 
\cite{bkpsv}.
The Euler equation in Sobolev spaces is examined in \cite{bl1}, where a norm-inflation result for the related vorticity equation provides the foundation for its ill-posedness.  
For the ill-posedness of the Burgers equation in $H^{3/2}$ we refer 
to Linares, Pilod and Saut \cite{lps}.
For more results on traveling wave solutions, well-posedness and 
other analytic and geometric properties of nonlinear evolution equations
we refer the reader to the following works and the references therein 
\cite{bc},
\cite{bl2},
\cite{cht},
\cite{cl},
\cite{cm},
\cite{dgh},
\cite{eep},
\cite{ely},
\cite{ey},
\cite{h},
\cite{hk},
\cite{hkm},
\cite{hmp},
\cite{kl},
\cite{kt},
\cite{lo},
\cite{mc},
\cite{mst},
\cite{ti},
\cite{w}.

\vskip0.1in
This  paper  is organized as follows.
In Section \ref{sec:2},  we construct the 
2-peakon solutions on the line having the properties 
described in Theorem \ref{NE-2peak-slns}.
We begin with the system of the four differential 
equations  defined by  \eqref{odes}  when $n=2$
and after making the change of the dependent 
variables $q=q_2-q_1$, $p=p_2-p_1$, 
$w=p_2+p_1$, and $z=p_1p_2$,  
we solve the resulting system and find explicit formulas for
$p, w$ and $z$ in terms of $q$
(see Proposition \ref{q-p-w-z-solutions}).
 For $q=q(t)$ we obtain a rather complicated
autonomous differential equation, which can
be dominated by a simpler one for which 
we can prove, by a comparison argument,
 that $q$ becomes zero (collision)
in finite time. Also, a precise estimate of 
the collision time is derived.
This is contained in  Proposition \ref{zero-of-q-lifespan}.
In Section \ref{sec:3}, we estimate the $H^s$-norm 
of the 2-peakon solutions constructed earlier
(see Proposition \ref{norm-estimate-NE}),
and in section \ref{sec:4} we choose 
the parameter appropriately so that 
both the lifespan (collision time) 
and the size of the 2-peakon solution 
at the initial time are small.
In Section \ref{sec:5}, we prove
norm-inflation and illposedness
for $5/4<s<3/2$. 
Then,  in Section \ref{sec:6},
we prove non-uniqueness for $s<5/4$
by showing that our 2-peakon solution
$u(t)$ converges in $H^s$ to an antipeakon $u(T)$,
which gives rise to a second solution for NE
having the same initial data.
Also, we explain the ill-posedness of NE for $s=5/4$. 
Finally, in Section \ref{sec:7} we prove our results on the circle.  We
use analogous arguments to those used on the line, with the necessary modifications
to account for the periodic environment.  A detailed outline of the periodic case can 
be found in Subsection \ref{subsec:1}.
%
%
%
%
%
%
%
%
%
%
%
\section{Construction of 2-peakon solutions}
\label{sec:2}
 \setcounter{equation}{0}
It can be shown (see \cite{hw2008}, \cite{gh})
 that the   2-peakon  
\begin{align}
\label{2-peak}
u(x,t) = p_1(t) e^{-|x-q_1(t)|}+p_2(t) e^{-|x-q_2(t)|},
\end{align}
is a solution of  NE if the positions $q_1$, $q_2$
and the momenta $p_1$, $p_2$ satisfy 
the following system of the four differential equations
\begin{equation}
\begin{split}
\label{odes-2peakon}
q'_1 & = \Big(p_1 +p_2 e^{-|q_1-q_2|}\Big)^2,
\\
q'_2 & = \Big(p_1 e^{-|q_2-q_1|}+p_2\Big)^2,
\\
p'_1 &= p_1p_2\Big(p_1+p_2 e^{-|q_1-q_2|}\Big)
        \cdot
\text{sgn}(q_1-q_2)  e^{-|q_1-q_2|},
\\
p'_2  &= p_1p_2\Big(p_1 e^{-|q_2-q_1|}+p_2 \Big)
      \cdot
      \text{sgn}(q_2-q_1)  e^{-|q_2-q_1|},
\end{split}
\end{equation}
where $\text{sgn}(x)$ is the standard sign function defined
to be 1 if $x>0$, $-1$ if $x<0$, and $0$ if $x=0$.
At this point we make our first observation. Since $q'_1\ge 0$ and 
$q'_2\ge 0$, both position are increasing with time.  Therefore
we cannot have the ``typical" peakon-antipeakon collision which is created from the peakon traveling in the positive direction and antipeakon traveling in the negative direction as
observed in the cases of the CH and DP equations.
Also,  we note that by translation we may assume that the 
initial positions $q_1$ and $q_2$ are symmetric, that is 
\begin{align}
\label{q-initial-data}
q_1(0)=-a
\quad
\text{and}
\quad 
q_2(0)=a,
\quad
\text{for some } a>0,
\end{align}
and,  at least for a while,  the difference of the positions is positive,
that is
\begin{align}
\label{q-diff-pos}
q(t)=q_2(t)-q_1(t)>0.
\end{align}
Thus, the last system takes the following simpler form
\begin{equation}
\begin{split}
\label{2peakon-DE-simple-k=2}
q_1' & = \Big(p_1 +p_2 e^{-q}\Big)^2, \\
q_2' & = \Big(p_1 e^{-q}+p_2\Big)^2, \\
p_1' &= -p_1p_2\Big(p_1+p_2 e^{-q}\Big) e^{-q},\\
p_2'  &= p_1p_2\Big(p_1 e^{-q}+p_2 \Big)e^{-q}.
\end{split}
\end{equation}
Furtheremore,  we shall assume  that at time $t=0$  the 
initial  momenta are 
\begin{align}
\label{b-choice}
 p_2(0)=b\gg 1, \qquad    p_1(0)= -(b+\delta),\quad \delta>0.
\end{align}
That is, the initial profile $u_0(x)=u(x,0)$ is the following asymetric antipeakon-peakon
\begin{align}
\label{antpeak-peak-data}
u_0(x) = -(b+\delta) e^{-|x+a|} + be^{-|x-a|},
\end{align}
which is displayed in the Figure 1.
\vskip0.2in
\hskip1in
\begin{minipage}{0.7\linewidth}
\hspace*{0cm}
\vspace*{0cm}
\begin{tikzpicture}[xscale=1,yscale=1]
%
%
\newcommand\X{7};
\newcommand\Y{2};
\newcommand\FX{11};
\newcommand\FY{11};
\newcommand\R{0.6};
\newcommand\B{2};
\newcommand\D{0.5};
\newcommand\A{1};
%
%
\draw[->,line width=1pt,black] (-5,0)--(5,0) 
node[above left] {\fontsize{\FX}{\FY}$x$};
\draw[->,line width=1pt,black] (0,-2.5)--(0,2.5) node[below right] {\fontsize{\FX}{\FY}$u_0$};
\draw[domain=-4:4, variable=\x, 
red, line width=1.5pt] 
plot ({\x},{-(\B+\D)*exp(-abs(\x+\A))+\B*exp(-abs(\x-\A))});
\draw[line width=1pt,black,dashed] 
({-(\A)},{-(\B+\D)+0.35}) 
node[] { } 
node[below,xshift=-1.7cm] {\fontsize{\FX}{\FY}
$p_1(0)=-(b+\delta)\approx u_0(-a)$}
--
({-(\A)},0) 
node[above, xshift=-.2cm] {\fontsize{\FX}{\FY}$q_1(0)=-a$} 
node[] {$\bullet$} ;
\draw[line width=1pt,black,dashed] 
({(\A)},{(\B+\D)-0.9}) 
node[] {}
node[above,xshift=1.2cm] {\fontsize{\FX}{\FY}
$u_0(a)\approx b=p_2(0)$}
--
({(\A)},0) 
node[below, xshift=.2cm] {\fontsize{\FX}{\FY}$q_2(0)=a$} node[] {$\bullet$};
\end{tikzpicture}
\end{minipage}
\vskip0.15in
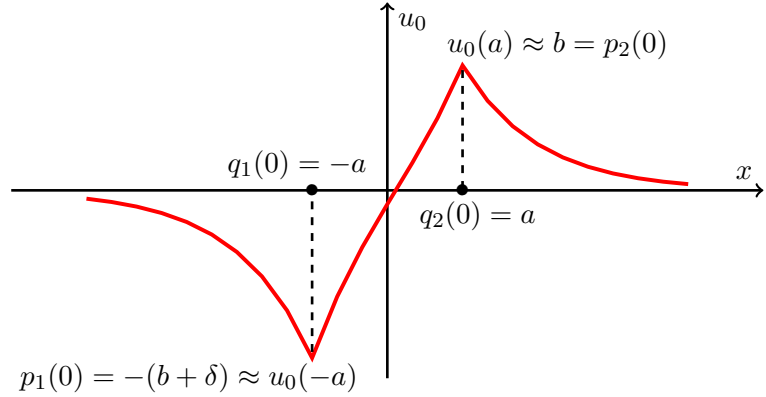
\centerline{Figure 1: Initial profile $u_0(x)$}
\vskip0.2in
Next, we shall solve the system of differential equations
\eqref{2peakon-DE-simple-k=2}
with initial data the antipeakon-peakon \eqref{antpeak-peak-data}
and prove that  there is a  collision in finite time.
To demonstrate this claim, it is more convenient  to work with the following
new dependent variables 
\begin{equation}
\begin{split}
\label{q-q-p-w-def}
q(t)&=q_2(t)-q_1(t), \quad \,\hspace{.02in} q(0)=2a>0,
\\
p(t)&=p_2(t)-p_1(t), \quad \, p(0)=2b+\delta>0,
\\
w(t)&=p_2(t)+p_1(t), \quad \hspace{-.01in} w(0)=-\delta<0,
\\
z(t)&=p_2(t)\cdot p_1(t), \quad \,\,\,\, z(0)=-b(b+\delta)<0.
\end{split}
\end{equation}
\vskip0.1in
\noindent
{\bf \large Deriving equations for $q$,
$p$, $w$ and $z$.}
Subtracting the first equation of the system 
\eqref{2peakon-DE-simple-k=2}
from the second we have 
\begin{align}
\notag
q'&=\Big(p_1 e^{-q}+p_2\Big)^2
-
 \Big(p_1 +p_2 e^{-q}\Big)^2,
 \\& =
    (p_2-p_1) (p_2 +p_1)(1-e^{-2q}),
    \notag
         \\& =
    p w(1-e^{-2q}).
    \label{deriv-of-q=q2-q1}
\end{align}
Next, we shall try to form differential equations  for $p$ and $w$
using the system \eqref{odes-2peakon}. Assuming $p_1<0$ and $p_2>0$,
at least for some time, for $p$ we have
\begin{align}
\label{p2-p1}
p'&= p_1p_2\Big(p_1 e^{-q}+p_2 \Big) e^{-q}
        +
        p_1p_2\Big(p_1+p_2 e^{-q}\Big)
 e^{-q},
\notag
  \\
  \notag&=  
 p_1p_2(p_2+p_1) e^{-q} (1+e^{-q}),
   \\&=  
zw e^{-q} (1+e^{-q}).
 \end{align}

For $w$ we have
\begin{align}
\notag
w'&=\Big(p_1 e^{-q}+p_2 \Big) e^{-q}
-
 p_1p_2\Big(p_1+p_2 e^{-q}\Big)  e^{-q},
\\&=
\notag
p_1p_2(p_2-p_1) e^{-q} (1-e^{-q}),
\\&=
zp e^{-q} (1-e^{-q}).
 \end{align}
\noindent
Finally,  for $z$ we have 
\begin{align}
\label{p2-times-p1}
z'&= p'_2\cdot p_1 + p_2\cdot p'_1,
\notag
\\&=  
p_1^2p_2\Big(p_1 e^{-q}+p_2 \Big) e^{-q}
       -
 p_1p_2^2\Big(p_1+p_2 e^{-q}\Big) e^{-q},
\notag
  \\\notag
  &=  
p_1p_2e^{-q}\Big[
 (p_1^2-p_2^2) e^{-q}
 \Big],
   \\&=  
-zwpe^{-2q}.
\end{align}
To summarize, we have the following system for $q$, $p$, $w$ and $z$
\begin{align}
\begin{split}
\label{q-de-NE}
q' &=  p w(1-e^{-2q}), \quad \,\,\,\,\hspace{.01in} q_0=q(0) \,\hspace{.01in} =2a>0,
\\  
p' &=   zw e^{-q} (1+e^{-q}),  \,\,\, p_0=p(0) \,\hspace{.01in}  =2b+\delta>0,
\\
 w'&= zp e^{-q} (1-e^{-q}), \,\,\,\hspace{.01in} w_0=w(0)=-\delta<0,
\\
 z'& =- zwpe^{-2q},  \quad \quad \,\,\,\,\,\hspace{.01in}  z_0=z(0)
 \,\hspace{.01in} =-b(b+\delta)<0.
\end{split}
\end{align}
In the following result we derive explicit formulas for $p$, $w$ and $z$
in terms of $q$. For $q$, we derive an autonomous differential 
equation,  which in turn, is dominated by a simpler such equation.
\begin{proposition}
[Solutions of transformed 2-peakon system]
\label{q-p-w-z-solutions}
The system of differential equations  \eqref{q-de-NE}
 has a unique smooth solution  $(q(t), p(t), w(t), z(t))$ in an  interval $[0, T)$,
for some $T>0$, such that  $z=z(t)$ is  decreasing
and in terms of $q$  is expressed  by the formula
\begin{align}                   
\label{z-q-form-NE}
z
=
\frac{-z_1}{\big(1-e^{-2q}\big)^{1/2}}<0,
 \quad
 \text{where}
  \quad
  z_1=b(b+\delta)\big(1-e^{-2q_0}\big)^{1/2},
\end{align}
$p=p(t)$ is  decreasing and as a function of $q$  is expressed  by the formula
\begin{align}
\label{p-q-form-NE}
p
=
\Big(
p_0^2 +
2z_1 \Big[
\frac{1+e^{-q}}{\sqrt{1-e^{-2q}}}
-
\frac{1+e^{-q_0}}{\sqrt{1-e^{-2q_0}}}
\Big]
\Big)^{1/2}>0,
\end{align}
and
$w=w(t)$ is  decreasing and as a function of $q$  is expressed  by the formula
\begin{align}
\label{w-q-form-NE}
w(t)
=
-\Big(
w_0^2 +
2z_1
 \Big[
\frac{\sqrt{1-e^{-2q_0}}}{1+e^{-q_0}}
-
\frac{\sqrt{1-e^{-2q}}}{1+e^{-q}}
\Big]
\Big)^{1/2}<0.
\end{align}
The difference of the positions $q=q(t)$ is decreasing and satisfies the initial value problem
\begin{align}
\notag
\label{q-ivp}
&q'=-f(q)
\doteq
-\Big(
w_0^2 +
2z_1
 \Big[
\frac{\sqrt{1-e^{-2q_0}}}{1+e^{-q_0}}
-
\frac{\sqrt{1-e^{-2q}}}{1+e^{-q}}
\Big]
\Big)^{\frac 12}\cdot
\\&
\hskip0.8in \;\;\;\;\;\;\;\;
\cdot
\Big(
p_0^2 +
2z_1 \Big[
\frac{1+e^{-q}}{\sqrt{1-e^{-2q}}}
-
\frac{1+e^{-q_0}}{\sqrt{1-e^{-2q_0}}}
\Big]
\Big)^{\frac 12}
\cdot
(1-e^{-2q}),
\\&
q(0)=q_0=2a>0.
\notag
\end{align}
Furthermore,  the initial value problem  \eqref{q-ivp} for $q$
is dominated by the simpler initial value problem 
\begin{align}
\label{q-simpler-dominant-ivp-NE}
q'=-g(q)
\doteq
-q_1 \big(1-e^{-2q}\big)^{3/4},
\quad
0<q(0)=2a<1/2,
\end{align}
where 
\begin{equation}
\label{q1-def-NE}
q_1
=
\delta 
\sqrt{2b(b+\delta)}\cdot q_0^{1/4}.
 \end{equation} 
\end{proposition}
{\bf  Proof.}
We begin by expressing $z$ in terms of $q$. Using the equation for $z'$ and $q'$,
we find 
\begin{align}
\notag
 \frac{z'}{q'}=\frac{- zpwe^{-2q}}{ p w(1-e^{-2q})},
 \quad
 \text{or}
  \quad
   \frac{z'}{z}=\frac{- e^{-2q}q'}{ (1-e^{-2q})}.
\end{align}
Since $z(0)<0$ we shall assume that $z(t)$ will remain negative.
Therefore, from the last relation we have 
\begin{align}
\notag
\frac{d}{dt}[\ln (-z)]
=
-\frac 12 \frac{d}{dt}[\ln (1-e^{-2q})].
\end{align}
Integrating from $0$ to $t$ gives
\begin{align}
\notag
\ln\Big[\frac{z(t)}{z_0}\Big]
=
-\frac 12 \ln \Big[\frac{1-e^{-2q}}{1-e^{-2q_0}}\Big].
\end{align}
Finally, solving for $z$ gives formula \eqref{z-q-form-NE},
which expresses $z$ in terms of $q$.

\vspace{.1in}

%
%
%
%
%
%
%
Next we express $p$ in terms of $q$. For this  we divide the 
equation for $p'$ by the equation for $q'$  and we get
\begin{align}
\notag
 \frac{p'}{q'}=\frac{zw e^{-q} (1+e^{-q})}{ p w(1-e^{-2q})},
 \quad
 \text{or}
  \quad
   pp'=z\cdot \frac{e^{-q} (1+e^{-q})q'}{ (1-e^{-2q})}.
\end{align}
Substituting into the above relation the formula for $z$
given by \eqref{z-q-form-NE}, we have
\begin{align}
\label{pp'-eq}
   pp'=\frac{-z_1}{\big(1-e^{-2q}\big)^{1/2}} \cdot \frac{e^{-q} (1+e^{-q})q'}{ (1-e^{-2q})}
   =
    \frac{-z_1(1+e^{-q})e^{-q} q'}{\big(1-e^{-2q}\big)^{3/2}}.
\end{align}
Furthermore, by making the change of variables,
$u=e^{-q(t)}$, we have $du=-e^{-q(t)}q'(t) dt$  and 
\begin{align*}
\int
    \frac{-(1+e^{-q})e^{-q} q'}{\big(1-e^{-2q}\big)^{3/2}}
dt
=
\int \frac{1+u}{\big(1-u^2\big)^{3/2}} du
=
\frac{1+u}{\big(1-u^2\big)^{1/2}}+C
=
\frac{1+e^{-q(t)}}{\big(1-e^{-2q(t)}\big)^{1/2}} +C.
\end{align*}
Therefore, relation \eqref{pp'-eq} reads as 
\begin{align}
\label{p-t-derivative-1}
\frac{d}{dt}\Big[\frac12 p^2\Big]
=
z_1 \frac{d}{dt}\Big[\frac{1+e^{-q}}{\sqrt{1-e^{-2q}}}\Big].
\end{align}
Integrating \eqref{p-t-derivative-1}  from $0$ to $t$ gives
\begin{align}
\notag
\frac12 \Big[p^2(t)-p_0^2\Big]
=
z_1 \Big[
\frac{1+e^{-q(t)}}{\sqrt{1-e^{-2q(t)}}}
-
\frac{1+e^{-q_0}}{\sqrt{1-e^{-2q_0}}}
\Big],
\end{align}
which, when solved  for $p$, gives formula \eqref{p-q-form-NE},
which expresses $p$ in terms of $q$.
\vskip0.1in

%
%
%
%
%
%
Finally, we express $w$ in terms of $q$. Dividing the 
equation for $w'$ by the equation for $q'$  gives 
\begin{align}
\notag
 \frac{w'}{q'}=\frac{zp e^{-q} (1-e^{-q})}{ p w(1-e^{-2q})},
 \quad
 \text{or}
  \quad
   ww'=z\cdot \frac{ e^{-q} (1-e^{-q})q'}{ (1-e^{-2q})}.
\end{align}
Now, substituting the formula for $z$
given by  \eqref{z-q-form-NE} into the above relation, we get
\begin{align}
\label{ww'-eq}
   ww'
   =
   \frac{-z_1}{\big(1-e^{-2q}\big)^{1/2}} \cdot \frac{e^{-q} (1-e^{-q})q'}{ (1-e^{-2q})}
      =
   \frac{-z_1(1-e^{-q})e^{-q} q'}{\big(1-e^{-2q}\big)^{3/2}}.
\end{align}
Furthermore,  making again the change of variables
$u=e^{-q(t)}$, we have
\begin{align*}
\int
   \frac{-(1-e^{-q})e^{-q} q'}{\big(1-e^{-2q}\big)^{3/2}}
dt
=
\int \frac{1-u}{\big(1-u^2\big)^{3/2}} du
=
-\frac{\sqrt{1-u^2}}{1+u}+C
=
-\frac{\sqrt{1-e^{-2q(t)}}}{1+e^{-q(t)}} +C.
\end{align*}
Therefore, relation \eqref{ww'-eq} reads as follows
\begin{align}
\label{intermediate-step-2}
\frac{d}{dt}\Big[\frac12 w^2\Big]
=
-z_1 \frac{d}{dt}\Big[
\frac{\sqrt{1-e^{-2q(t)}}}{1+e^{-q(t)}}
\Big].
\end{align}
Integrating \eqref{intermediate-step-2} from $0$ to $t$ gives
\begin{align}
\notag
\frac12 \Big[w^2(t)-w_0^2\Big]
=
z_1 \Big[
\frac{\sqrt{1-e^{-2q_0}}}{1+e^{-q_0}}
-
\frac{\sqrt{1-e^{-2q(t)}}}{1+e^{-q(t)}}
\Big].
\end{align}
Solving for $w$  while  taking into consideration
that $w(t)<0$ in the choice of sign, 
gives formula \eqref{w-q-form-NE},
which expresses $w$ in terms of $q$.

\vspace{.1in}

Concerning  the differential equation for $q$, we begin from its equation
$q'= wp(1-e^{-2q})$ and substituting for $w$ and
$p$ their expressions \eqref{w-q-form-NE} and \eqref{p-q-form-NE},
we obtain the desired autonomous
 initial value problem \eqref{w-q-form-NE}.
Next, we observe that 
\begin{align}
\label{w-simplification}
\frac{\sqrt{1-e^{-2q_0}}}{1+e^{-q_0}}
-
\frac{\sqrt{1-e^{-2q}}}{1+e^{-q}}
\ge 0, 
\quad
0\le q\le q_0,
\end{align}
and also  that
\begin{align}
\label{p-simplification}
p_0^2  -
2z_1 \frac{1+e^{-q_0}}{\sqrt{1-e^{-2q_0}}}
\ge 0
\iff
\frac{(2b+\delta)^2}{ 2b(b+\delta)}\ge 1+e^{-q_0}.
\end{align}
In fact, condition \eqref{p-simplification}  is implied by the stronger
condition
\begin{align}
\notag
\frac{(2b+\delta)^2}{ 2b(b+\delta)}\ge 2
\iff
4b^2+4b\delta+\delta^2>4b^2+4b\delta
\iff
\delta^2>0,
\quad
\text{which is true.}
\end{align}
Now, using  \eqref{w-simplification} and  \eqref{p-simplification} 
we see that  the function $f(q)$ in the right-hand side  of the differential equation \eqref{q-ivp} can be bounded from below 
by
\begin{align}
\notag
f(q)
&\ge
\Big( 
w_0^2
\Big)^{\frac 12}
\cdot
\Big(
2z_1 \Big[
\frac{1+e^{-q}}{\sqrt{1-e^{-2q}}}
\Big]
\Big)^{\frac 12}
\cdot
(1-e^{-2q})
\\&
\notag
= 
\delta
\cdot
\Big(
2b(b+\delta) \big(1-e^{-2q_0}\big)^{1/2}\Big)^{\frac 12}
\Big(
 \Big[
\frac{1+e^{-q}}{\sqrt{1-e^{-2q}}}
\Big]
\Big)^{\frac 12}
\cdot
(1-e^{-2q}).
\end{align}
Using the bounds  $1+e^{-q}>1$ and  
$1-e^{-2q_0}\ge q_0$, for $0\le q_0\le 1/2$, which follow from 
the following simple but useful approximation 
\begin{align}
\notag
\frac x2 \le 1-e^{-x}\le x \iff 1-e^{-x} \simeq x, 
\quad \text{if} 
\quad  0\le x\le 1,
  \end{align}
  we have 
\begin{align}
\notag
f(q)
\ge
\delta
\sqrt{2b(b+\delta)} \cdot q_0^{1/4}
\cdot
\big(1-e^{-2q}\big)^{3/4}
\doteq
g(q).
\end{align}
Therefore, defining 
$
q_1
\doteq
\delta
\sqrt{2b(b+\delta)}\cdot q_0^{1/4}
$
we see that the complicated initial value problem  for $q$
given in
\eqref{q-ivp} is dominated by the 
simpler one shown in \eqref{q-simpler-dominant-ivp-NE}.
\qed

\vskip0.1in
Next we move our attention to the study of the solution $q(t)$
 of the initial value problem stated in Proposition \ref{q-p-w-z-solutions}. From the formulas for $p$ and $w$, we see that
  they blow-up at a zero of $q$. Therefore, 
 the lifespan of  our  2-peakon solution is equal to
the first such  zero. 
 The following result, which is applicable
 to the simpler dominant initial value problem 
 \eqref{q-simpler-dominant-ivp-NE}
 proves existence of  a zero and  provides an estimate for it size
 in terms of the initial data.

\begin{proposition}
[Zero of $q$]
\label{zero-of-q-lifespan}
If $r<1$ then for given $q_0\in (0, 1/2)$ and $q_1>0$
the solution to the initial value problem
\begin{align}
\label{q-ivp-r}
\frac{dq}{dt}  
=
-g_r(q)
\doteq
-q_1
\left(1-e^{-2q}\right)^r,
\quad
q(0)=q_0,
\end{align}
which begins positive and is decreasing, becomes zero in finite time
$T$ given by 
\begin{align}
\label{T-lifespan-r}
T
=
\int_0^{q_0}
 \frac{dq}{g_r(q)}
 =
\frac{1}{q_1}
\int_0^{q_0}
 \frac{dq}{\left(1-e^{-2q}\right)^r}
\simeq
 \frac{1}{1-r}\, \frac{q_0^{1-r}}{q_1}.
\end{align}
\end{proposition}
A key ingredient in proving Proposition \ref{zero-of-q-lifespan}
is the following elementary result that 
compares solutions of the initial value problem
\eqref{q-ivp-r} for different values of $r$.
It states  that a bigger $r$ correspond to 
a bigger solution.
\begin{lemma}
[Comparison principle]
\label{qr-comparison}
 If  $r_1$  and $r_2$ are two values of $r$ 
 such that   $r_1\le r_2$, then the corresponding solutions 
 $q_{r_1}(t)$ and $q_{r_2}(t)$ to the initial value problem \eqref{q-ivp-r}
 with the same initial  data $q_0$  satisfy $q_{r_1}(t) \le  q_{r_2}(t)$. That is,
\begin{align}
\notag
 r_1\le r_2
  \quad 
  \Longrightarrow
  \quad 
 q_{r_1}(t) \le  q_{r_2}(t).
\end{align}
\end{lemma}
\noindent
{\it Proof.}  It follows from the fact that $r_1\le r_2$ implies
\begin{align}
\notag
-q_1\left(1-e^{-2q}\right)^{r_1} \le - q_1\left(1-e^{-2q}\right)^{r_2}. \qed
\end{align}
\medno
{\bf Remark.}  We note that for $r\ge 1$ the solution
to the initial value problem \eqref{q-ivp-r} has no zero.
In fact, for $r=1$ it  reads as follows
\begin{align}
\notag
\frac{dq}{dt} 
 =-q_1
 \left(1-e^{-2q}\right),
  \quad 
 q(0)=q_0.
\end{align}
Integrating this equation, gives the explicit formula
\begin{align}
\notag
q(t)=\frac 12 \ln\Big[1+(e^{2q_0} -1)e^{-q_1\, t}\Big]
\doteq q_1(t).
\end{align}
From this formula we see that the solution $q(t)$ exists for all $t\ge 0$,
is positive for all times and decreases to zero as $t$ goes to $\infty$.
Thus when $r=1$ then $q_1(t)$ has no zero in finite time.
Since, by the comparison principle the solution 
$q_r(t)$ that corresponds to an $r>1$ is greater to $q_1(t)$,
we conclude that $q_r(t)$ has no zero in finite time if $r>1$.
Therefore, the lifespan $T$ is equal to $\infty$ if $r\ge 1$.

\vskip0.1in
\noindent
{\bf Proof of Proposition \ref{zero-of-q-lifespan}.}
We begin with the case $r\le 0$.
  When  $r=0$ then
 our initial value problem  \eqref{q-ivp}  become the following 
 simple one
$
q'(t)
 =-q_1
 $,
$
 q(0)=q_0,
$
whose solution is
\begin{align}
\notag
q(t)= q_0 -q_1\, t\doteq q_0(t),
\end{align}
which has a zero at $T=q_0/q_1$. Thus, 
by the comparison Lemma \ref{qr-comparison} the solution 
$q_r(t)$ that corresponds to an $r<0$  is smaller to $q_0(t)$,
and therefore has a zero in finite time. In fact, it is smaller than
$q_0/q_1$. This proves existence of zero for $q_r(t)$ when $r\le 0$.

 \vskip0.1in
 \noindent
{\bf  Existence of a zero for $q(t)$ if $0<r<1$:} 
To prove existence of zero of $ q_{r}(t)$ for $0<r<1$,
it suffices to  do so
under  the additional condition
\begin{align}
\label{r-not-(n-1)/n}
r \ne \frac{n-1}{n}, \quad \text{for all} \quad n=1, 2, 3, \cdots.
\end{align}
In fact,  if $r$ were of the form $\frac{n-1}{n}$ then we could choose
another $r_2 \in (0, 1)$  which is not of this form  and 
$r<r_2$. Then,  by the comparison Lemma \ref{qr-comparison},
proving the  existence of  a zero for $q_{r_2}(t)$ implies 
existence of  a zero for $q_{r}(t)$.  So,  from now on we 
shall assume that $r$ satisfies condition 
\eqref{r-not-(n-1)/n}.
Therefore, there is a positive integer $n\ge 2$  such that 
\begin{align}
\label{r-double-ineq}
\frac{n-2}{n-1}<r<\frac{n-1}{n}.
\end{align}
It turns out that for proving existence 
of a zero of $q=q_{r}(t)$,
we need its  $n$-th order Taylor polynomial approximation  at $t=0$.
Differentiating  equation \eqref{q-ivp} $n$ times, we arrive at the formula
\begin{align}
\label{q-nth-deriv}
q^{(n)}(t) 
&=
q_1^n c_n(r)
\left(1-e^{-2q(t)}\right)^{nr-(n-1)}
+
 q_1^n\,\sum_{j=1}^{n-1}c_j(r) \left(1-e^{-2q(t)}\right)^{nr-(j-1)},
\end{align}
where 
\begin{align}
\label{q-c-n}
c_n(r)
=
(-1)^n 2^{n-1} r(2r-1)\cdots \big([n-1]r-[n-2]\big)
\end{align}
and $c_j(r)$ for $j=1, \cdots, n-1$ are coefficients depending on $r$.
Also, we obtain the following formula for 
the $(n+1)$-th derivative of $q$
\begin{align}
\label{q-(n+1)th-deriv}
\hskip-0.07in
q^{(n+1)}(t) 
=
q_1^{n+1} \hskip-0.02in c_{n+1}(r) \hskip-0.02in  \left(1-e^{-2q(t)}\right)^{(n+1)r-n}
\hskip-0.08in
+
 q_1^{n+1}\sum_{j=1}^{n} \hskip-0.01in c_j(r)
 \hskip-0.02in \left(1-e^{-2q(t)}\right)^{(n+1)r-(j-1)},
\end{align}
where 
\begin{align}
\label{q-c-n+1}
c_{n+1}(r)
=
(-1)^{n+1} 2^{n} r(2r-1)\cdots \big(nr-[n-1]\big), 
\end{align}
and again $c_j(r)$ for $j=1, \cdots, n$ are coefficients depending on $r$.  Therefore, the $n$-th order Taylor polynomial approximation 
of $q(t)$ at $t=0$ is given by
\begin{align}
\label{nth-taylor-approx}
q(t)
=
q_0 + q'(0)t+\frac{q''(0)}{2!} t^2 +\frac{q^{(3)}(0)}{3!} t^3+\cdots+
\frac{q^{(n)}(0)}{n!} t^n
+
\frac{q^{(n+1)}(\tau)}{(n+1)!} t^{n+1},
\end{align}
where $0\le \tau \le t$. Next, we shall show
that the coeficients $c_n(r)$ and $c_{n+1}(r)$
defined by \eqref{q-c-n} and \eqref{q-c-n+1}
have the same sign, which
is the key ingredient 
for proving the existence of 
a zero for $q(t)$.
We prove this 
claim by considering the two cases
possible, $n$ even and $n$ odd.
We begin with the case of $n$ even.
In this case, using the first part
of inequality  \eqref{r-double-ineq} that $r$ satisfies, 
we see that $(n-1)r>n-2$
and this implies that $c_n(r)$ is a positive number. Also,
 using the second part of inequality   \eqref{r-double-ineq} 
 we see that $nr<n-1$, which implies that $c_{n+1}(r)>0$
 is a positive number too. Furthermore,
 in the expression of $q^{(n)}(0)$ the first term  
$ q_1^n c_n(r)
\left(1-e^{-2q_0}\right)^{nr-(n-1)}$ is the dominant term 
for $q_0$ small enough since the exponent $nr-(n-1)$
is negative while the exponents of $(1-e^{-2q_0})$ appearing 
in all other terms of
the sum  \eqref{q-nth-deriv} are positive. Thus we can conclude
\begin{align}
\label{q-nth-deriv-positive}
q^{(n)}(0)>0,
\quad
\text{if $n$ is even and  $q_0$ is small enough.}
\end{align}
Similarly,  in the expression of $q^{(n+1)}(t)$, 
the first term  
$ q_1^{n+1} c_{n+1}(r)
\left(1-e^{-2q(t)}\right)^{(n+1)r-n}$ is the dominant term 
for $q_0$ small enough, since the exponent $(n+1)r-n$
is negative while all the exponents of $(1-e^{-2q(t)})$ appearing in the sum  \eqref{q-(n+1)th-deriv} are positive,
except the one that corresponds to $j=n$ which 
has exponent $(n+1)r-(n-1)$ which is negative.
However,  $\left(1-e^{-2q(t)}\right)^{(n+1)r-(n-1)}$
 is  dominated by  $\left(1-e^{-2q(t)}\right)^{(n+1)r-n}$,
 for $q(t)\le q_0$ small enough.
 Thus, we also have
\begin{align}
\label{q-(n+1)-deriv-positive}
q^{(n+1)}(\tau)>0,
\quad
\text{if $n$ is even and  $q_0$ is small enough.}
\end{align}
In the case that $n\ge 2$ is an odd positive integer
then the signs change due to the fact $(-1)^n=-1$
and using the same reasoning as in the even case 
we  obtain that 
\begin{align}
\label{q-nth-deriv-negative}
q^{(n)}(0)<0,
\quad
\text{if $n$ is odd and  $q_0$ is small enough,}
\end{align}
and 
\begin{align}
\label{q-(n+1)-deriv-negative}
q^{(n+1)}(\tau)<0,
\quad
\text{if $n$ is odd and  $q_0$ is small enough.}
\end{align}
Now, we are ready to prove the existence of zero for $q(t)$.
First we consider the case that  $n$ is an odd number.
 Then, using the $n$-th order Taylor polynomial approximation 
 \eqref{nth-taylor-approx}
and the conditions \eqref{q-nth-deriv-negative},
\eqref{q-(n+1)-deriv-negative} we obtain that 
\begin{align}
\notag
q(t)
\le
q_0 + q'(0)t+\frac{q''(0)}{2!} t^2 +\frac{q^{(3)}(0)}{3!} t^3+\cdots+
\frac{q^{(n)}(0)}{n!} t^n,
\quad
\text{for all}
\quad
 t\ge 0.
\end{align}
Furthermore, since for large $t$ the term $\frac{q^{(n)}(0)}{n!} t^n$ 
dominates and $q^{(n)}(0)<0$ we have that the $n$-th order Taylor polynomial approximation of $q(t)$ will become negative, thus
crossing the $t$-axis. This forces $q(t)$ to have a zero
at some positive time $T$, which is the desired
conclusion.

Finally, we prove the existence of zero for $q(t)$ in the even case. This is done by contradiction. In fact, if $q(t)>0$ for all $t>0$ 
then our differential equation 
$q'(t)  = -q_1  \left(1-e^{-2q(t)}\right)^r$ 
implies that $q(t)$ is decreasing for all $t>0$
and therefore 
\begin{align}
\label{q-less-than-q0}
q(t)
\le
q_0 
\quad
\text{for all}
\quad
 t\ge 0.
\end{align}
However, if $n$ is even, then using the
 $n$-th order Taylor polynomial approximation 
of $q(t)$ at $t=0$, which is  given by \eqref{nth-taylor-approx},
and  conditions \eqref{q-nth-deriv-positive} and
\eqref{q-(n+1)-deriv-positive}, we have that 
\begin{align}
\label{q-greater-than-poly}
q(t)
\ge
q_0 + q'(0)t+\frac{q''(0)}{2!} t^2 +\frac{q^{(3)}(0)}{3!} t^3+\cdots+
\frac{q^{(n)}(0)}{n!} t^n,
\quad
\text{for all}
\quad
 t\ge 0.
\end{align}
Inequality \eqref{q-greater-than-poly} leads to
a contradiction
because for large $t$ the  term  $\frac{q^{(n)}(0)}{n!} t^n$
dominates the Taylor polynomial approximation. Thus, 
there is some large time  $T>0$ such that 
\begin{align}
\label{q-greater-than poly-T}
q(t)
\ge
q_0  +\frac 12 \frac{q^{(n)}(0)}{n!} T^n>2q_0,\quad 
0\le t\le T,
\end{align}
which contradicts inequality \eqref{q-less-than-q0}.
This argument completes the proof of the existence of zero for $q(t)$
when $0<r<1$.

%
%
%
%
%
%
%
%
%
%
%
%
\vskip0.1in
\noindent
{\bf Estimating  the  zero $T$ of the position $q(t)$ when $r<1$.}
Let $T$ be the zero of the solution $q(t)$ of our initial value 
problem \eqref{q-ivp}, which is: 
$q'(t)
 =
 -q_1
 \left(1-e^{-2q(t)}\right)^r,
 \,
q(0)=q_0.
$
Integrating it from $0$ to $T$ we have
\begin{align}
\notag
\int_0^T
 \frac{q'(t)}{\left(1-e^{-2q(t)}\right)^r} dt
 =
-q_1T.
\end{align}
Then,  making the substitution $q=q(t)$, and
using the initial and terminal conditions $q(0)=q_0$
and $q(T)=0$,  we obtain the following formula for $T$
\begin{align}
\label{T-formula-1} 
T
=
\frac{1}{q_1}
\int_0^{q_0}
 \frac{dq}{\left(1-e^{-2q}\right)^r}
 T
\simeq
\frac{1}{q_1}
\int_0^{q_0}
 \frac{dq}{q^r}
 =
 \frac{1}{1-r}\, \frac{q_0^{1-r}}{q_1}.
\end{align}
Above, we used the estimate 
 $q\le 1- e^{-2q}\le 2q$, if $0\le q\le \frac{1}{2}$.
This completes the proof of Proposition \ref{zero-of-q-lifespan}.
 \qed

\vskip0.1in
Applying Proposition \ref{zero-of-q-lifespan}  with 
$r=3/4$  we obtain the following result  for the
the zero the initial value problem   \eqref {q-ivp}
and the lifespan of our 2-peakon solution $u$.
\begin{corollary}
[Zero of $q$ and lifespan of $u$]
\label{zero-of-q-lifespan-of-u}
If $0<q_0<1/2$ and $b$,  $\delta$ satisfy condition 
\eqref{b-choice} 
 then the solution  to the initial value problem
 \eqref {q-ivp}
begins positive, is decreasing, and becomes zero in finite time
$T$ given by
 {\small
\begin{align}
\label{T-lifespan-3/4}
T
=
\int_0^{q_0}
 \frac{dq}{f(q)}
 \le
 \frac{1}{q_1}
\int_0^{q_0}
 \frac{dq}{\left(1-e^{-2q}\right)^{3/4}}
\simeq
 \frac{q_0^{1/4}}{q_1}
 \simeq
  \frac{q_0^{1/4}}{
\delta
\sqrt{2b(b+\delta)}\cdot q_0^{1/4}
  }
   \simeq
  \frac{1}{
\delta
\sqrt{2b(b+\delta)}
  }.
  \end{align}
  }
\end{corollary}
\noindent
{\bf Proof.} The existence and uniqueness  of  the 
solution follows 
from the fundamental ODE theorem since $f(q)$ is 
a smooth function. That $q(t)$ is decreasing follows 
from the fact that $q'=-f(q)<0$. Finally, that $q(t)$ becomes zero
in finite time follows from the fact that our initial
value problem   \eqref {q-ivp} is dominated 
by the initial value problem \eqref{q-simpler-dominant-ivp-NE}
for which Proposition 
\ref{zero-of-q-lifespan}
is applicable
with $r=3/4$. Therefore, estimate \eqref{T-lifespan-r}
gives \eqref{T-lifespan-3/4}, and this completes 
the proof of the lemma. \qed
\vskip0.1in

%

%
%
%
%
%
%
%
The  properties of our special 2-peakon solutions
are summarized in the following Theorem
and are a consequence of
Proposition \ref{q-p-w-z-solutions} and 
Corollary  \ref{zero-of-q-lifespan-of-u}.
\begin{theorem}
[Construction of 2-peakon solutions]
\label{cor:NE-2-peakon}
For given $0<a\le 1/4$ and  $b$,  $\delta$  satisfying  condition \eqref{b-choice}  the initial value problem
for the positions $q_1, q_2$ and the momenta $p_1, p_2,$
\begin{equation}
\label{NE-2peakon-ivp}
\begin{split}
q_1' & = \Big(p_1 +p_2 e^{-q}\Big)^2, \quad\quad\quad\quad \hspace{.01in} q_1(0)=-a,
\\
q_2' & = \Big(p_1 e^{-q}+p_2\Big)^2, \quad\quad\quad\quad q_{2} (0)=a>0,
\\
p_1' &= -p_1p_2\Big(p_1+p_2 e^{-q}\Big) e^{-q},
 \,\, p_1(0)=-b-\delta,
\\
p_2'  &= p_1p_2\Big(p_1 e^{-q}+p_2 \Big)e^{-q}, \quad \hspace{.01in} p_{2}(0)=b,
\end{split}
\end{equation}
\end{theorem}
%
%
%
%
%
\begin{minipage}{.5\textwidth} %
{\it
has a unique smooth solution $(q_1, q_2, p_1, p_2)(t)$
with a finite lifespan  $T$, which is the zero of $q=q_2-q_1$, and which satisfies the estimate 
\begin{equation}
\label{NE-life-T-est}
T
   \lesssim
  \frac{1}{
\delta
\sqrt{2b(b+\delta)}
}.
\end{equation}
Furthermore, we have
\begin{align}
\notag
&p_1=\frac{w-p}{2}<0, \,\, decreasing, \,\,  \\ \notag
&\lim_{t \to T^{-}} p_1(t) = -\infty, \,\, 
  \text{and}  \,\,  -p_1\simeq p\simeq q^{-1/4},
\end{align}
and 
\begin{align}
\notag
&p_2=\frac{w+p}{2}>0, \,\, increasing, \,\,  \\ \notag
 &\lim_{t \to T^{-}} p_2(t) = \infty, \,\, 
  \text{and}  \,\,  p_2\simeq p  \simeq q^{-1/4},
\end{align}
where $p$ and  $w$ are given in Proposition \ref{q-p-w-z-solutions}.
Also,  $w=p_1+p_2$ is decreasing from 
$ w_0 < 0$ to $w_T$, where $w_T
\doteq
 \lim_{t \to T^{-}} w(t)$, that is
\begin{align}
\label{wT-def}
w_T
 =
-\Big(
\delta^2 +
2   b(b+\delta)
(1-e^{-2a})
\Big)^{\frac 12}.
\end{align}
Finally,  the  2-peakon  
\begin{align}
\notag
u(x,t) = p_1(t) e^{-|x-q_1(t)|}+p_2(t) e^{-|x-q_2(t)|},
\end{align}
is NE solution for $x\in\rr$, $0<t<T$,
with the following asymmetric  antipeakon-peakon initial profile
\begin{align}
\notag
u(x,0) = -(b+\delta) e^{-|x+a|}+b e^{-|x-a|}, \,
x\in\rr.
\end{align}
}
\end{minipage} 
%
%
%
%
\hskip0.2in
\begin{minipage}{.5\textwidth} %
\begin{tikzpicture}[xscale=1,yscale=0.5]
%
%
\newcommand\X{7};
\newcommand\Y{2};
\newcommand\FX{11};
\newcommand\FY{11};
\newcommand\R{0.6};
\newcommand\B{2};
\newcommand\D{0.5};
\newcommand\A{1};
%
%
\draw[->,line width=1pt,black] (-1,0)--(3.75,0) 
node[above left] {\fontsize{\FX}{\FY}$t$};
\draw[->,line width=1pt,black] (0,-11)--(0,11);
\draw[domain=0:2, range=-5:5,variable=\x, smooth, red, line width=1.5pt,->] 
plot ({\x},{exp(\x)+\B+1});
\draw[domain=0:2, range=-5:5,variable=\x, smooth, red, line width=1.5pt,->] 
plot ({\x},{-(exp(\x)+\B+1)-\D});
\draw[domain=0:2, range=-5:5,variable=\x, smooth, red, line width=1.5pt] 
plot ({\x},{2*\A*exp(-\x)-.19});
\draw[domain=0:2, range=-5:5,variable=\x, smooth, red, line width=1.5pt] 
plot ({\x},{2*\A*exp(-\x)-3.19});
\draw[dashed,line width=1pt]
(0,{\B+2}) 
node[] {$\bullet$}
node[left]  {\fontsize{\FX}{\FY}$p_2(0)=b$}
(0,{-\B-2-\D}) 
node[] {$\bullet$}
node[xshift=-1.4cm]  {\fontsize{\FX}{\FY}$p_1(0)=-b-\delta$}
(0,{2*\A-.19}) 
node[] {$\bullet$}
node[left]  {\fontsize{\FX}{\FY}$2a$}
(0,{2*\A-3.19}) 
node[] {$\bullet$}
node[left]  {\fontsize{\FX}{\FY}$w_0=-\delta$}
(0,{2*\A-5}) 
node[] {$\bullet$}
node[left]  {\fontsize{\FX}{\FY}$w_T$}
--
(2.3,{2*\A-5})
(2.15,11) 
--
(2.15,0)
node[] {$\bullet$}
node[below right] {\fontsize{\FX}{\FY}$T$}
--
(2.15,-11);
\node[rotate=55] at (1.2,5) {\fontsize{\FX}{\FY}$p_2(t)
\approx q^{-\frac{1}{4}}$};
\node[rotate=-55] at (1.2,-5.8) {\fontsize{\FX}{\FY}$p_1(t)\approx -q^{-\frac{1}{4}}$};
\node[] at (1.2,1.7) {\fontsize{\FX}{\FY}$q(t)$};
\node[] at (1.2,-1.7) {\fontsize{\FX}{\FY}$w(t)$};
\end{tikzpicture}
\vskip.15in
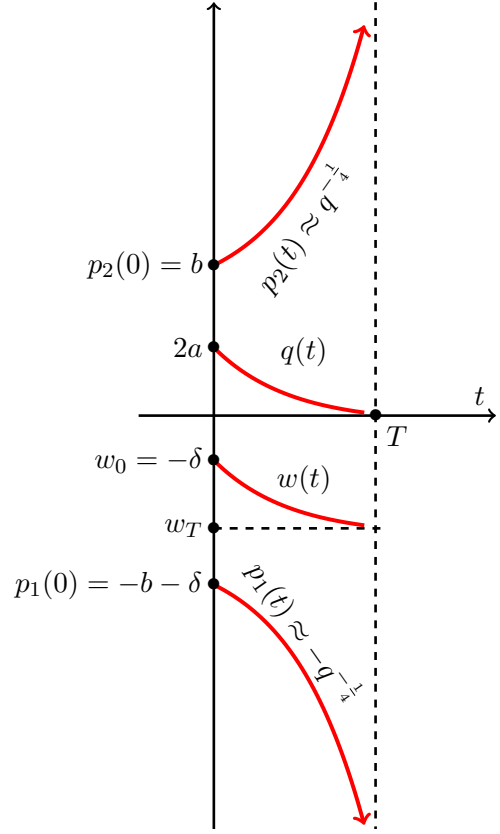
\centerline{Figure 2: Graphs of $
p_1$, $p_2$, $q$ and 
$w$ }
\end{minipage}
%
%

%
%
%
%
%
%
\section{Calculating the Norm}
 \label{sec:3}
\setcounter{equation}{0}

\begin{proposition}
\label{norm-estimate-NE}
Let $u(t)$ be the two-peakon solution 
to the NE equation. Then on $[0,T)$ we have 
\begin{align}
\label{u-norm-NE}
\|u(t)\|_{H^s}^2
=
16 r(t) p_1^2(t)     Q_s(q)
+
4c_s\big(1-r(t)\big)^2 p_1^2(t),
\quad
\text{with}
\quad 
r(t)\doteq -\frac{p_2(t)}{p_1(t)},
\end{align}
where    $c_s=  \int_\rr  (1+\xi^2)^{s-2}  d\xi$  and  
$Q_s(q)$, which  is given below,
satisfies  the estimates:
\begin{align}
\label{Qs-estimate}
 Q_s(q)
 \doteq
 \int_\rr
(1+\xi^2)^{s-2} \sin^2 \bigg(\frac{q\xi}{2} \bigg) d\xi
\simeq
\begin{cases}
q^{3-2s},     \quad \;\quad\quad 1/2 < s < 3/2, \\
q^2 \cdot \ln(1/q),   \;\;\;\, s = 1/2, \\
q^2,\quad \quad \quad \quad \quad s < 1/2.
\end{cases}
\end{align}
\end{proposition}

{\bf Proof.} Since $\widehat{e^{-|x|}}(\xi)=2/(1+\xi^2)$ we have that 
the Fourier transform of 
\begin{align}
\notag
u(x,t) = p_1 e^{-|x-q_1|} + p_2 e^{-|x-q_2|}
\end{align}
is given by
\begin{align}
\notag
\widehat{u}(\xi,t)
&=
\frac{2 p_1 e^{-i \xi q_1}}{1+\xi^2} 
+
\frac{2p_2 e^{-i\xi q_2}}{1+\xi^2}
=
\frac{2}{1+\xi^2}  \cdot 
 p_1 e^{-i \xi  q_1} \cdot \Big(1 + \frac{p_2}{p_1} e^{-i \xi  q} \Big).
\end{align}
Taking the square of the $H^s$ norm of this quantity and factoring out $p_1^2$,
we obtain
\begin{align}
\label{norm-of-u-1}
\|u(t)\|_{H^s}^2
&=
4 p_1^2 
\int_\rr
(1+\xi^2)^{s-2} \Big| 1 + \frac{p_2}{p_1} e^{-i\xi q} \Big|^2 d\xi.
\end{align}
Using Proposition \ref{q-p-w-z-solutions} we see that
\begin{align}
\notag
r=r(t)\doteq -\frac{p_2(t)}{p_1(t)}=\frac{p+w}{p-w}<1,
\quad
\text{and}
\quad
 r(t) \nearrow 1 \text{ as } t \nearrow T.
\end{align}
Next, using  $r$ we write \eqref{norm-of-u-1} as follows
\begin{align}
\label{u-norm-squared}
\|u(t)\|_{H^s}^2
=
4 p_1^2 
\int_\rr
(1+\xi^2)^{s-2} \Big| 1-r e^{-i\xi q}  \Big|^2 d\xi.
\end{align}
Expanding out the square under the integral
in \eqref{u-norm-squared}, we have
\begin{align}
|re^{i q \xi} - 1|^2
\label{formula-square-exp}
 &=
(1-r)^2+ 4 r \sin^2 \bigg(\frac{q\xi}{2} \bigg).
\end{align}
Substituting \eqref{formula-square-exp} into \eqref{u-norm-squared}
\begin{align}
\notag
\|u(t)\|_{H^s}^2
=
16 r p_1^2 
\int_\rr
(1+\xi^2)^{s-2} \sin^2 \bigg(\frac{q\xi}{2} \bigg) d\xi
+
4(1-r)^2 p_1^2 
\int_\rr
(1+\xi^2)^{s-2}  d\xi,
\end{align}
or
\begin{align}
\notag
\|u(t)\|_{H^s}^2
=
16 r p_1^2     Q_s(q)
+
4c_s(1-r)^2 p_1^2,
\end{align}
where
\begin{align}
\label{cs-def}
c_s=  \int_\rr  (1+\xi^2)^{s-2}  d\xi
\quad
\text{and}
\quad
Q_s(q)
=
\int_\rr
(1+\xi^2)^{s-2} \sin^2 \bigg(\frac{q\xi}{2} \bigg) d\xi.
\end{align}
%
%
%
%
Now, we see that  to prove Proposition \ref{norm-estimate-NE}
 it suffices to show that  for $0<q< 1/8$
 we have the following estimate
\begin{align}
\notag
 Q_s(q)
\simeq
\begin{cases}
q^{3-2s},  \quad \quad \quad \quad \;\quad\quad 1/2 < s < 3/2, \\
q^2 \cdot \ln(1/q), \quad\quad \quad  \;\;\; s = 1/2, \\
q^2,\quad \quad \quad \quad  \quad \quad \quad \quad s < 1/2.
\end{cases}
\end{align}
Starting with the integrand for $Q_s$ from \eqref{cs-def} and 
making the  change of variables 
$x = q\xi$, which gives  $dx = qd\xi$, 
we can write  $Q_s(q)$ as  
\begin{align}
Q_s(q)
&=
2q^{3-2s} \int_0^\infty (q^2 +  x^2)^{s-2} \sin^2(x/2) dx
=
2q^{3-2s} \big[  I_1 + I_2\big],
\label{Hs-I1-I2}
\end{align}
where
\begin{align}
\notag
I_1 \doteq \int_0^1 \frac{ x^2}{(q^2 +  x^2)^{2-s}} dx
\quad
\text{and}
\quad
I_2 \doteq \int_1^\infty \frac{ \sin^2(x/2)}{(q^2 +  x^2)^{2-s}} dx.
\end{align}
If $ s<3/2$ then  the integral  $I_2$ is bounded since
\begin{align}
I_2
=
\int_{1}^\infty  \frac{\sin^2(x/2) }{(q^2 +  x^2)^{2-s}} dx 
\lesssim
\int_{1}^\infty 
x^{2s - 4} dx 
=\frac{1}{3-2s}.
\label{I2-up-bd}
\end{align}
Also,  when $s>1/2$ we have the following upper bound for $I_1$
\begin{align}
\label{I1-up-bd}
I_1
\le
\int_{0}^1  x^{2s-2}dx 
=\frac{1}{2s-1}.
\end{align}
Furthermore, for any $s<2$  we have 
\begin{align}
I_1
\ge
2^{s-2}
\int_0^{1} x^2 dx 
= 2^{s-2} \cdot \frac 13.
\label{I1-low-s-below-2}
\end{align}
Combining  \eqref{I1-low-s-below-2},  
 \eqref{I2-up-bd},
 and \eqref{Hs-I1-I2} gives 
\begin{align}
\notag
\|f_q\|_{H^s} 
\simeq
q^{\frac{3}{2}-s}, 
\quad
\text{if}
\quad 
  1/2 < s < 3/2.
\end{align}
{\bf The case $s<1/2$:}  Since $I_1+I_2$ is bounded below by $I_1$
and $0<y<1$  we have 
\begin{align}
I_1+I_2
\gtrsim
\int_0^{q} \frac{x^2}{(q^2 +  x^2)^{2-s}}dx 
\gtrsim
q^{2s-4}  \int_0^{q} x^2 dx
\simeq
q^{2s - 1}.
\label{I-low-bd}
\end{align}
Combining  \eqref{I-low-bd} and  \eqref{Hs-I1-I2} gives
\begin{align}
 Q_s(q) \gtrsim 
q^{3-2s} \cdot q^{2s - 1}=q^2,
\quad
\text{if}
\quad 
  s < 1/2.
\label{b1/2-I-low-bd}
\end{align}
To prove  the reverse of  inequality \eqref{b1/2-I-low-bd} 
we obtain an upper bound  for $I_1$.
For this argument, we let $z=x/y$ and get
\begin{align*}
I_1 \doteq \int_0^1 \frac{ x^2}{(q^2 +  x^2)^{2-s}} dx
=
q^{2s-1}
\int_0^{1/q} \frac{ z^2}{(1 +  z^2)^{2-s}} dz
\le
q^{2s-1}
\int_0^{\infty} \frac{ 1}{(1 +  z^2)^{1-s}} dz.
\end{align*}
Since the last integral converges if $2(1-s)>1$,
which is equivalent to $s<1/2$, we see that
it is equal to a finite constant $c_s$.
\noindent
Combining this  fact together with  \eqref{Hs-I1-I2} and
  \eqref{I2-up-bd}   we  have
\begin{align}
\notag
 Q_s(q)
\lesssim
q^{3-2s}[ q^{2s-1} +1]
\lesssim
q^2,
\quad
\text{if}
\quad 
  s<1/2,
\end{align}
which together with \eqref{I-low-bd} gives
\begin{align}
\notag
 Q_s(q)
\simeq
y^2,
\quad
\text{if}
\quad 
  s<1/2.
\end{align}
%


%
{\bf The case $s=1/2$:}  We observe that
\begin{align}
\notag
I_1 
&=
\int_0^1 \frac{x^2}{(q^2 + x^2)^{3/2}} dx
=
\ln \Big( \sqrt{q^2 + 1} + 1 \Big)- \frac{1}{\sqrt{q^2 + 1}}
+
\ln (1/q).
\end{align}
{\bf Upper Bound.} 
From here we begin by removing the middle 
term and using the fact that $y < 1/4$ in the first term.  We get
\begin{align}
I_1 
\le 
\ln(\sqrt{2} + 1) + \ln(1/q).
\label{1/2-I1-up}
\end{align}
Substituting \eqref{1/2-I1-up} back into \eqref{Hs-I1-I2} 
and taking into account estimate \eqref{I2-up-bd} for $I_2$ we have
\begin{align}
Q_{\frac 12}(s)
&\lesssim
q^{2} \ln(1/q).
\label{1/2-I1-up-bd}
\end{align}
{\bf Lower Bound.}
Using the fact that $2 \ln(1/q) > 1$ we have
\begin{align}
\notag
I_1 
\ge
\ln(1/q) + (\ln(2) - 1) [2 \ln(1/q)]  
=
(2\ln(2) - 1) \ln(1/q).
\end{align}
We therefore arrive at
\begin{align}
\label{1/2-I1-low-bd}
Q_{\frac 12}(s)
\gtrsim
q^2 \cdot I_1  
\gtrsim
q^2 \ln(1/q). 
\end{align}
Putting these upper and lower bounds \eqref{1/2-I1-up-bd} and 
\eqref{1/2-I1-low-bd}  together and taking the square root of both sides
of the equation gives the desired result of
\begin{align}
\notag
Q_{\frac 12}(s) \simeq y^2 \cdot \ln(1/q).   \qed
\end{align}
%

%
%
%
%
%
%
%
\section{ Small lifespan and initial data}
 \label{sec:4}
\setcounter{equation}{0}
We begin by assuming that 
\begin{align}
\notag
 p_2(0)=b\gg1 \text{ and }-p_1(0)=  b+\delta,\, \delta>0,
\end{align}
so that  the conditions for the existence
of our 2-peakon with the lifespan estimate
\eqref{T-lifespan-3/4}
hold. Then,  we have the following.

\medno
{\bf Lifespan Estimate.}   For given $\ee>0$,
we need to find $b>1$ such that $T<\ee$. 
Since, by Proposition \ref{q-p-w-z-solutions} we have
\begin{align}
\notag
T
\lesssim
\frac{1}{
\delta
\sqrt{2b(b+\delta)}}
\le
\frac{1}{
\delta b},
\end{align}
we must have
\begin{align}
\label{T-small-1}
\frac{1}{
\delta b}
\le \ee \iff  b\ge \delta^{-1} \ee^{-1}.
\end{align}

\medno
{\bf Initial Data Estimate.} Now, for the same $\ee>0$
we need to find  $q_0<1/8$ such that
 $\|u_0\|_{H^s} <\ee$. 
 For this 
argument
 we use Proposition \ref{norm-estimate-NE},
 from which we have
\begin{align}
\notag
\|u(0)\|_{H^s}^2
&=
16 r(0) p_1^2(0)     Q_s(q_0)
+
4c_s\big(1-r(0)\big)^2 p_1^2(0),  \quad  r(t)\doteq -\frac{p_2(t)}{p_1(t)},
\\&\notag
=
16 b(b+\delta)    Q_s(q_0)
+
4c_s\delta^2,
\end{align}
which in turn gives
\begin{align}
\notag
\|u(0)\|_{H^s}^2
\le
32 b^2    Q_s(q_0)  
+
4c_s\delta^2.
\end{align}
%
%
%
%
{\bf Case $1/2 < s < 3/2$:} Then  by Proposition \ref{norm-estimate-NE}
we have  $Q_s(q_0) \lesssim q_0^{3-2s}$ and therefore
\begin{align}
\notag
\|u(0)\|_{H^s}^2
\le
C_sb^2   q_0^{3-2s}
+
4c_s\delta^2.
\end{align}

To demonstrate   $\|u_0\|_{H^s} <\ee$, it suffices 
to choose $q_0$ and $\delta$ such that
$
C_sb^2   q_0^{3-2s}
+
4c_s\delta^2
\le \ee^2$, or
\begin{align}
\notag
4c_s\delta^2  \le \frac{\ee^2}{2}
\quad
\text{ and }
\quad
C_sb^2   q_0^{3-2s} \le \frac{\ee^2}{2}.
\end{align}
The first  inequality holds if
\begin{align}
\label{delta-choice}
\delta \le \frac{\ee}{2\sqrt{2c_s}}.
\end{align}
Taking into consideration \eqref{delta-choice} and 
\eqref{T-small-1}, the second  inequality holds if
\begin{align}
\notag
q_0^{3-2s} \le \frac{\ee^2}{2C_sb^2 }
\le
\frac{\ee^2}{2C_s \delta^{-2} \ee^{-2} }
=
\frac{\delta^{2}\ee^4}{2C_s}
\le
\frac{\ee^{2}\ee^4}{8c_s\cdot 2C_s},
\end{align}
or 
\begin{align}
\notag
q_0 
\le
\Big(
\frac{\ee^6}{16c_sC_s}
\Big)^\frac{1}{3-2s}.
\end{align}
%
%
%
%
%
%
%
%
{\bf Case $s \le 1/2$:} For such a Sobolev exponent $s$
we have   $\|u(0)\|_{H^s} \le \|u(0)\|_{H^1}$. 
This combined with  Proposition \ref{norm-estimate-NE},
which tells us that   $Q_1(q_0) \lesssim q_0$,
gives
\begin{align}
\notag
\|u(0)\|_{H^s}^2
\le
\|u(0)\|_{H^1}^2
\le
C_1b^2   q_0
+
4c_1\delta^2.
\end{align}
Thus  $\|u_0\|_{H^s} <\ee$  if 
$q_0$ and $\delta$ satisfy the inequalities 
\begin{align}
\notag
4c_1\delta^2  \le \frac{\ee^2}{2}
\quad
\text{ and }
\quad
C_1b^2   q_0 \le \frac{\ee^2}{2}.
\end{align}
These inequality holds if 
\begin{align}
\notag
\delta \le \frac{\ee}{2\sqrt{2c_1}}
\quad
\text{and}
\quad
q_0 
\le
\frac{\ee^6}{16c_sC_s}.     
\end{align}
%

%
%

\vskip0.2in
\section{Norm-Inflation and illposedness for $5/4<s<3/2$}
 \label{sec:5}
\setcounter{equation}{0}
From Proposition \ref{norm-estimate-NE}
we have
\begin{align}
\label{u-norm-NE-fin}
\|u(t)\|_{H^s}^2
=
16 r(t) p_1^2(t)     Q_s(q)
+
4c_sp_1^2(t)\big(1-r(t)\big)^2,
\end{align}
where the estimate for $Q_s$ is given in \eqref{Qs-estimate}.  Also, using
Theorem \ref{cor:NE-2-peakon}  we  have
\begin{equation}
\notag
p_1^2(t)\simeq q^{-1/2}(t),
\quad
\text{and}
\quad 
p_2^2(t)\simeq q^{-1/2}(t),
\quad
\text{for $t$ close to $T$}.
\end{equation}
Next, we see that
\begin{align}
\notag
r=r(t)\doteq \frac{p_2(t)}{-p_1(t)}
\simeq 
\frac{q^{-1/4}}{q^{-1/4}}
\simeq 
1,
\quad
 \text{ as } t \nearrow T,
\end{align}
and
\begin{equation}
\notag
p_1(t)\big(1-r(t)\big)
=
p_1(t)\Big(1+\frac{p_2(t)}{p_1(t)}\Big)
=
p_2(t)+p_1(t)=w(t).
\end{equation}
Also, we have
\begin{equation}
\label{w-lim-fin}
\lim_{t\to T} p_1^2(t)\big(1-r(t)\big)^2
=
\lim_{t\to T} w^2(t)
=
\delta^2 +
2b(b+\delta)
\cdot
(1-e^{-q_0}).
\end{equation}
Therefore, the first term of 
\eqref{u-norm-NE-fin} can be estimated by
\begin{align}
\notag
16 r(t)p_1^2(t)     Q_s(q)
\simeq
\begin{cases}
q^{\frac 52-2s},  & 1/2 < s < 3/2, \\
q^{\frac 32} \cdot \ln(1/q), & s = 1/2, \\
q^{\frac 32}, & s < 1/2.
\end{cases}
\end{align}
Combining the last estimate with the fact  $\frac 52-2s=0 \iff s=\frac 54$
we see that
\begin{align}
\label{Q-limts}
\lim_{t\to T}  16r(t)p_1^2(t)    Q_s(q)
=
\begin{cases}
\infty\quad
\text{(inflation)},     &    5/4 < s < 3/2, \\
\text{may not exist},  &  s = 5/4, \\
0,  &  s < 5/4.
\end{cases}
\end{align}
Finally, using the limits \eqref{Q-limts} and \eqref{w-lim-fin}
from formula \eqref{u-norm-NE-fin} we conclude that
\begin{equation}
\label{limit-H^s-norm-u}
\lim_{t\to T} \|u(t)\|_{H^s}^2
=
\begin{cases}
\infty\quad
\text{(inflation)},     &    5/4 < s < 3/2, \\
\text{may not exist},  &  s = 5/4, \\
4c_s\Big[
\delta^2 +
2b(b+\delta)
\cdot
(1-e^{-q_0})\Big],  &  s < 5/4.
\end{cases}
\end{equation}
Therefore when   $5/4<s <3/2$ we have norm inflation
and ill-posedness for the Novikov equation. \qed

%
%
%
%
%
%
%
\section{Non-Uniqueness for $s<5/4$}
 \label{sec:6}
 \setcounter{equation}{0}

In this section, we prove that once we take the Sobolev exponent to be less that $5/4$, the Novikov
equation admits non-unique solutions.  
\begin{theorem}[Non-uniqueness] 
For $s < 5/4$  NE  admits non-unique solutions.
\label{NE-non-uniqueness}
\end{theorem}
Our proof of non-uniqueness revolves around examining the behavior of the limit as 
$t \to T^-$ of the 2-peakon
solution $u$ with initial data given in \eqref{antpeak-peak-data}.  Once we take 
the Sobolev exponent to be $s < 5/4$, this limit exists, and it is a single antipeakon.  
The non-uniqueness
then can be realized by taking a single antipeakon traveling wave that which at time $T$ has the same
profile as $\lim_{t \to T^-} u(x,t)$.  From this point, a change of variables can recast this scenario as two solutions arising from the same initial data.  To proceed with this argument, we 
begin by examining the pointwise limit, then the $L^r$ limit, and 
finally use these results in addition to the generalized Dominated Convergence Theorem to establish
the $H^s$ limit.

\begin{proposition} [Pointwise limit]
\label{u-pointwise-limit} 
For each $x \in \rr$ we have 
\begin{align}
\label{sol-t-limit}
\lim_{t \to T^-} u(x,t) = w_T e^{-|x - q_T|} \doteq v_T(x),
\end{align}
where $w_t$ is given by \eqref{wT-def},
that is 
$w_T =-\Big(
\delta^2 +
2   b(b+\delta)
(1-e^{-2a})
\Big)^{\frac 12}<0$, and
\begin{align}
\notag
& q_T \doteq \lim_{t \to T^-} q_1(t) = \lim_{t \to T^-}  q_2(t).
\end{align}
\end{proposition}
For proving Proposition \ref{u-pointwise-limit} we shall 
need the following elementary result.
\begin{lemma}
\label{lem-p-pointwise-limits}
Given our functions $p_1, p_2$ and $q=q_1-q_1$
the following limits hold as $t \to T^-$:
\begin{align}
\label{p-pointwise-limits}
 \lim_{t \to T^-}    p_j(t) (1 - e^{-q(t)}) = 0, \quad j =1,2 . 
\end{align}
\end{lemma}
{\bf Proof.}  Using the estimates
\begin{align}
\notag
p_1^2 \simeq q^{-1/2} \quad \text{and} \quad p_2^2 \simeq q^{-1/2},
\end{align}
and the  inequality
$1 - e^{-x} < x$ for $x \in [0,1]$
we have
\begin{align*}
\lim_{t \to T^-}   |p_j(t) (1 - e^{-q(t)})| 
\le
\lim_{t \to T^-}   |p_j(t)| \cdot |q(t)|  
\lesssim
\lim_{t \to T^-}  |q^{-1/4} (t) | \cdot |q(t)| = 0. \qed
\end{align*}
\vskip0.1in
\noindent
{\bf Proof of Proposition \ref{u-pointwise-limit}.}  As we are working with a pointwise limit, we consider the cases 
$x \ge q_T$ and  $x < q_T$ separately so that we can evaluate
the absolute values $|x-q_j|$ in the definition of the 2-peakon
solution $u$.

\vspace{.1in}

{\bf Case $x \ge q_T$.} Since $q_1\le q_2\le q_T$,
we have $x-q_j\ge 0$ and therefore 
\begin{align}
\notag
u(x,t) 
&=
p_1(t) e^{-|x - q_1(t)|} + p_2(t) e^{- |x - q_2(t)|}
=
e^{-x} \cdot  \Big( p_1(t) e^{q_1(t)} + p_2(t) e^{q_2(t)} \Big).
\end{align}
Next, we will rewrite $u$ in such a way so as to utilize Lemma \ref{lem-p-pointwise-limits}.  
We have
\begin{align}
\notag
u(x,t) 
&=
e^{-x} \cdot 
e^{q_2(t)}
\cdot
\Big(
- p_1(t)  (1 - e^{-q(t)}) +  w(t)  \Big).
\end{align}
Finally  taking the limit as $t \to T^-$ of $u$ and 
using \eqref{p-pointwise-limits} we get
\begin{align}
\label{u-limit-x-above-qT}
\lim_{t \to T^-} u(x,t) = w_T \cdot e^{-x + q_T},
\quad x\ge q_T.
\end{align}
\vskip0.05in
{\bf Case $x < q_T$.}  We follow essentially the same strategy as in the previous case, simply correcting for signs.  Since $x$ is fixed and  $q_1 \le q_2 \le q_T$,
we see  that  after some 
time $t_0$ we must have $x < q_1(t) \le q_2(t)  \le q_T$.  
Therefore, for $t>t_0$ we have $x-q_j<0$  and $u$ can be written as
\begin{align}
\notag
u(x,t) 
&=
e^x \cdot 
e^{-q_1(t)} 
\Big( w(t) - p_2(t) (1 - e^{-q(t)})  \Big). 
\end{align}
Thus taking the limit as $t \to T^-$ of $u$
and using again Lemma  \ref{lem-p-pointwise-limits}
we obtain
\begin{align}
\label{u-limit-x-below-qT}
\lim_{t \to T^-} u(x,t) 
=
w_T \cdot e^{x - q_T},
\quad
 \quad x<q_T.
\end{align}
Combining \eqref{u-limit-x-above-qT} and \eqref{u-limit-x-below-qT} 
we conclude that the 2-peakon solution $u(t)$ has a limit 
as $t\to T^-$, which is given by the antipeakon  \eqref{sol-t-limit}.
\qed

\vspace{.1in}
We next examine the the limit of $u$ 
in $L^r$ topology.   
\begin{proposition}[Convergence in $L^r$] 
\label{Lr-conv}
For our antipeakon-peakon solution $u$ to NE, we have  
\begin{align}
\notag
\lim_{t \to T^-} \| u(x,t) - v_T(x) \|_{L^r} = 0,
\quad 
\text{for}
\quad
1 \le r < 4.
\end{align}
\end{proposition}

{\bf Proof.}
As we will need to evaluate the absolute values in the exponents, 
we note that the order of the peaks positions  of $u(x,t)$
and $v_T(x)$ is $q_1(t) < q_2(t) < q_T$. 
 We now expand the $L^r$ norm as
\begin{align}
\notag
\| u(x,t) - v_T(x) \|_{L^r}^r
\doteq
 I_1(t) +  I_2(t) +  I_3(t) +  I_4(t).
\end{align}
where the integrals $ I_j(t)$  have their domains determined by $q_1 < q_2 <q_T$, that is
\begin{align*}
& I_1(t) \doteq 
\int_{-\infty}^{q_1(t)}
|u(x,t) - v_T(x)|^r dx, 
&I_3(t) \doteq 
\int_{q_2(t)}^{q_T}
|u(x,t) - v_T(x)|^r dx,   
\\
& I_2(t) \doteq 
\int_{q_1(t)}^{q_2(t)}
|u(x,t) - v_T(x)|^r dx,  
&I_4(t) \doteq 
\int_{q_T}^{\infty}
|u(x,t) - v_T(x)|^r  dx.
\end{align*}
{\bf Evaluating $I_1$.} 
Calculating the integral, we have
\begin{align}
\notag
I_1(t) 
&=
\frac{ e^{rq_1(t)}}{r} \cdot 
\Big|
p_1 e^{- q_1 (t) } + p_2 e^{-q_2 (t) }  - w_T e^{ - q_T} \Big|^r .
\end{align}
In order to proceed with evaluating the limit, we observe the following identity
\begin{align}
\notag
p_1 e^{- q_1 (t) } + p_2 e^{-q_2 (t) }  - w_T e^{ - q_T}
&=
e^{-q_2}  p_1 (e^{q}  - 1) + w(t) e^{-q_2 (t) }  - w_T e^{ - q_T}.
\end{align}
We can now evaluate the limit as
\begin{align}
\notag
\lim_{t \to T^-}
I_1(t) 
&=
\lim_{t \to T^-}
\frac{ e^{r q_1(t)}}{r}  \cdot 
\Big|
e^{-q_2}  p_1 (e^{q}  - 1) + [w(t) e^{-q_2 (t) }  - w_T e^{ - q_T}]\Big|^r = 0.
\end{align}
{\bf Evaluating $I_2$.} 
Using the  Jensen's inequality
$
|a_1 +\cdots +a_n|^r \le n^r(|a_1|^r +\cdots + |a_n|^r)
$
together with   $e^{-|x-q_j (t)|}\le 1$, \, $e^{-|x-q_T(t)|}\le 1$
and   $|p_j|\simeq q^{-1/4}$, we have 
\begin{align*}
\lim_{t \to T^-}
I_2(t)
&\lesssim
\lim_{t \to T^-}
(q^{1- \frac{r}{4}}(t) + q(t)|w_T|^r)
= 0 \quad \text{(assuming $r < 4$)}.
\end{align*}
{\bf Evaluating $I_3$.}  After evaluating the absolute values inside the exponential,
and using the identity
\begin{align*}
p_1 e^{-x + q_1 } + p_2  e^{-x+q_2 }  - w_T e^{x - q_T}
&=
e^{-x}  e^{q_2}  p_1  ( e^{-q  }  - 1 )  + e^{-x} w e^{q_2} 
- e^{x} w_T e^{ - q_T},
\end{align*}
an application of Jensen's inequality gives us
\begin{align}
\notag
I_3(t) 
&\lesssim
\int_{q_2(t)}^{q_T}
\Big|e^{-x}  e^{q_2}  p_1 (t) ( e^{-q (t) }  - 1 )   \Big|^r dx
+
\int_{q_2(t)}^{q_T}
\Big|
e^{-x} w(t) e^{q_2(t)} - e^{x} w_T e^{ - q_T}
\Big|^r
dx.
\end{align}
We see that for the first term in this sum, we have
\begin{align}
\notag
\int_{q_2(t)}^{q_T}
\Big|e^{-x}  e^{q_2(t)}  p_1 (t) ( e^{-q (t) }  - 1 )   \Big|^r dx, 
&\lesssim
(q_T - q_2(t) )  
\cdot \Big|e^{q_2 (t)}  p_1 (t) ( e^{-q (t) }  - 1 )   \Big|^r.
\end{align}
For the second term of this sum, we use  the 
fact that $|e^{-x} w e^{q_2} - e^{x} w_T e^{ - q_T}|^r \lesssim 1$
and H\"{o}lder's inequality to get
\begin{align}
\notag
\int_{q_2(t)}^{q_T}
\Big|
e^{-x} w(t) e^{q_2(t)} - e^{x} w_T e^{ - q_T}
\Big|^r
dx
&\lesssim
(q_T - q_2 (t)). 
\end{align}
Putting these estimates together, we can now evaluate the limit of $I_3$ as $t \to T$ via
\begin{align}
\notag
\lim_{t \to T} I_3
\lesssim
\lim_{t \to T}
\Big(
(q_T - q_2(t) )  
\cdot \Big|e^{q_2 (t)}  p_1 (t) ( e^{-q (t) }  - 1 )   \Big|^r +
(q_T - q_2(t) \Big)
=
0.
\end{align}

{\bf Evaluating $I_4$.}    This term is handled in precisely the same fashion as $I_1$.
Performing the integration gives us
\begin{align}
\notag
I_4(t) 
&=
\frac{e^{-rq_T}}{r}
\cdot
\Big|
p_1 e^{ q_1 (t)} + p_2 e^{ q_2 (t)}  - w_T e^{ q_T} \Big|^r.  
\end{align}
Rewriting the expression inside 
of the absolute value gives us
\begin{align}
\notag
p_1 e^{ q_1 (t)} + p_2 e^{ q_2 (t)}  - w_T e^{ q_T} 
=
e^{q_2(t)} p_1 (e^{-q(t)} - 1) + w(t) e^{q_2(t)} - w_T e^{ q_T}.
\end{align}
Therefore, using the above identity along with the triangle inequality yields
\begin{align}
\notag
\lim_{t \to T^-}
I_4(t) 
&\le
 \lim_{t \to T^-} 
\frac{ e^{- r q_T}  }{r}
\cdot
\Big( |e^{q_2(t)} p_1 (e^{-q(t)} - 1)  |
+  |w(t) e^{q_2(t)} - w_T e^{ q_T} |\Big)^r = 0.
\end{align}
{\bf  Summarizing the $L^r$ convergence, $1 \le r < 4$.}  
As we have  computed $\lim_{t \to T^-} I_j (t) = 0$ for $j = 1,2,3,4$ it immediately follows
that
\begin{align}
\notag
\lim_{t \to T^-} \| u(x,t) - v_T(x) \|_{L^r}^r
&=
\lim_{t \to T^-} \Big(
I_1(t) + I_2(t) + I_3(t) + I_4(t) \Big) = 0. \qed
\end{align}
\begin{corollary}
\label{u-L2-conv}
As $t$ goes to $T$ our  2-peakon solution $u(t)$
 converges in  $H^s$, $s\le 0$,  to the antipeakon 
 $v_T=w_Te^{-|x-q_T|}$.
\end{corollary}

\vskip0.05in
Now that we have successfully established pointwise and $L^r$ convergence, we are ready
to move on to a much stronger result that is of interest in and of itself. 
 As $t \to T$, the antipeakon-peakon solution converges
to a single solitary antipeakon in $H^s$, for $s < 5/4$.  

\begin{theorem}
[Convergence in $H^s$]
\label{nu-hs} 
 For $s < 5/4$, our 2-peakon solution $u(t)$ 
converges to the antipeakon  $v_T$ in $H^s$, i.e.  
\begin{align}
\lim_{t \to T} \| u(t) - v_T\|_{H^s} = 0.
\end{align}
\end{theorem}

{\bf Proof.}
We will begin by simplifying the $H^s$ norm of $u(x,t) - v_T(x)$.
We have
\begin{align}
\notag
\|u(t) - v_T \|_{H^s}^2
&=
\int_\rr
(1 + \xi^2)^s
| \widehat{u}(\xi, t) - \widehat{v_T}(\xi)|^2
d\xi.
\end{align}
Our objective, when taking the limit,  will be to move the limit inside of the integral.  Thus, 
the first thing we should verify is whether pointwise, we have
\begin{align}
\notag
(1 + \xi^2)^s
| \widehat{u}(\xi, t) - \widehat{v_T}(\xi)|^2 \to 0 \quad \text{as} \quad t \to T^-.
\end{align}
We have
\begin{align}
\notag
| \widehat{u}(\xi, t) - \widehat{u}(\xi, T)|
&\le
\| u(x,t) - u(x,T) \|_{L^1}.
\end{align}
As we have proved that $u(x, t) \to u(x,T)$ in $L^1$, we have
\begin{align}
\notag
\lim_{t \to T^-}
(1 + \xi^2)^s
| \widehat{u}(\xi, t) - \widehat{v_T}(\xi)|^2 
&=
0.
\end{align}
Next, we will define the bounding
functions that will allow us to apply the generalized Dominated Convergence Theorem (gDCT).  We set
\begin{align}
\notag
f_t (\xi) \doteq
(1 + \xi^2)^s
|\widehat{u}(\xi, t) - \widehat{v_T}(\xi)|^2
&\le
4 (1+\xi^2)^s 
\Big(
|\widehat{u}(\xi, t)|^2 + |\widehat{v_T}(\xi)|^2  \Big) \doteq g_t (\xi). 
\end{align}
Next, we need to establish that the $g_t$'s, have a pointwise limit $g$.  The most
obvious candidate for $g$ is
\begin{align}
\notag
g(\xi)
=
8(1 + \xi^2)^s |\widehat{v}_T(\xi)^2|
\end{align}
Indeed, using the laws of limits, we have
$\widehat{u}(\xi,t) \to  \widehat{v}_T(\xi) $ pointwise in $\xi$
implies
$g_t \to g$ pointwise in $\xi$.
To finish satisfying the hypotheses of 
gDCT, we must now
establish the integral properties of the $g_t$'s.
We have
\begin{align}
\notag
\lim_{t \to T}
\int_\rr g_t (\xi) d \xi
&=
\lim_{t \to T}
4
\int_\rr  
(1+\xi^2)^s 
|\widehat{u}(\xi, t)|^2
d \xi
+
\lim_{t \to T}
4
\int_\rr (1+\xi^2)^s 
|\widehat{v_T}(\xi)|^2 d\xi 
\\
\notag
&
\overset{\eqref{limit-H^s-norm-u}}{=}
32 c_s w_T^2.
\end{align}
where the left limit uses the  hypothesis that $s<5/4$.
Furthermore, $g$ is  integrable and  
 \begin{align}
 \notag
 \int_\rr g(\xi) d\xi =
 \int_\rr 8 |v_T(\xi)^2| d\xi =
  8 \cdot \|v_T\|_{H^s}^2 = 8 \cdot 4 c_s w_T^2 = 32 c_s w_T^2.
 \end{align}
Therefore, we have
$
 \int g_t\to \int g.
 $
We now see that the hypotheses for  gDCT  are satisfied.
Thus, we can conclude  
\begin{align}
\notag
\lim_{t \to T}
\int_\rr f_t (\xi) d\xi
=
\int_\rr f(x) d\xi,
\end{align}
which written more explicitly tells us that
\begin{align}
\notag
\lim_{t \to T} 
\int_\rr 
(1 + \xi^2)^s
|\widehat{u}(\xi, t) - \widehat{v_T}(\xi)|^2 d\xi
&=
0.
\end{align}
Thus, we can conclude that as $t \to T$, we have $u(t) \to v_T$ in $H^s$.
\qed

\vspace{.1in}
{\bf Proof of Theorem \ref{NE-non-uniqueness}.} 
Translating the NE  1-peakon solution \eqref{NE-peak}
by $x_0$ and choosing the minus sign we obtain the following antipeakon solution for NE
\begin{equation}
\notag
v(x,t)= - \sqrt{c}\,  e^{-|(x-x_0)-ct|}, \quad \text{for any} \quad c>0
\quad
\text{and} \quad  x_0\in \rr.
\end{equation}
Choosing
\begin{equation}
\notag
c= w_T^2, 
\quad \text{and} \quad
x_0 = q_T - w^2_T T 
\end{equation}
we obtain the NE antipeakon solution 
\begin{equation}
\notag
v(x,t)= - \sqrt{w_T^2}\,  e^{-|(x- q_T + w^2_T T )-w^2_T\,t|}.
\end{equation}
Since at $t=T$ we have 
\begin{equation}
\notag
v(x,T)= w_T e^{-|x- q_T|}=u(x, T),
\end{equation}
we see that we have constructed two different NE solutions,
which belong in $H^s$, $s<5/4$, and agree at $t=T$.  
From here, a  change of 
variables can recast these two solutions as initial value problems at time $t=0$.
This proves failure of uniqueness in this range
of Sobolev spaces.  \qed

\vskip0.1in
{\bf The case $s=5/4$:}  If $s=5/4$ then
there are two possibilities. Either our 2-peakon 
solution $u(t)$ does not converge in $H^{5/4}$ in which case
we can prove (by a standard argument) that continuity
of the solution map fails, or  $u(t)$ converges in $H^{5/4}$
and has limit $u(T)$ (since this is the limit for lower
Sobolev exponents). In the second case, 
we have non-uniqueness like in Theorem 
\ref{NE-non-uniqueness}.
This result completes the proof of both of Theorems
\ref{NE-2peak-slns} and \ref{NE-ILP}
in the non-periodic case.

%
%
%
%
%
%
%
\section{The Periodic Case}
 \setcounter{equation}{0}
 \label{sec:7}
\subsection{Outline of the proofs in the periodic case} 
\label{subsec:1}
The proofs of Theorems \ref{NE-2peak-slns} and \ref{NE-ILP}
have been demonstated  on the line and we now  
present these proofs on the circle, $\tor = \rr/2\pi\zz$. The key ingredient is using a periodic version 
of the peakon. 
In Subsection \ref{subsec:2},  we construct the 
2-peakon solutions on the circle having the properties 
described in Theorem \ref{NE-2peak-slns}.
In Subsection \ref{subsec:3}, we estimate the $H^s$-norm 
of the 2-peakon solutions 
and in  Subsection \ref{subsec:4} we choose 
the parameters  so that 
both the lifespan 
and the size of the 2-peakon solution 
at the initial time are simultaneously small.
In Subsection \ref{subsec:5}, we prove
norm-inflation and illposedness
for $5/4<s<3/2$. 
Finally, in Subsection \ref{subsec:6},
we prove non-uniqueness for $s<5/4$
and 
explain the ill-posedness of NE for $s=5/4$. 
%

%
%
%
%
%
%
%
%
\subsection{Construction of 2-peakon solutions on the circle} 
\label{subsec:2}
The 2-peakon solutions to the periodic version of Novikov's equation are similar to those on the real-line with the caveat that the peak is generated by periodizing the hyperbolic cosine rather that using the exponential of the negative absolute value.  
The following equations are taken from \cite{gh} and 
\cite{hm}
and can also be derived in a straightforward fashion.

\vspace{.1in}
The periodic Novikov 2-peakon solutions  are of the form
\begin{align}
\label{per-NE-2-peakon}
u(x,t) = p_1(t) \cosh([x-q_1(t)]_p - \pi)  + p_2(t) \cosh([x-q_2(t)]_p - \pi),
\end{align}
where $[ \cdot ]_p$  periodizes our function and is defined  by the floor
\begin{align}
\label{x-per-def}
[x]_p = x -  2\pi \Big \lfloor \frac{x}{2\pi}  \Big\rfloor.
\end{align}
We see that $u$ solves NE if the momenta $p_1, p_2$ and the positions $q_1, q_2$
satisfy the following system of ODEs, which  can be obtained by using Theorem 1.2 of \cite{hm}
with the choice of parameters  $a = 0$, $b=3$.  The $4 \times 4$ system we get is
\begin{align}
\label{per-system-orig}
\begin{split}
&q_1' = p_1^2 [ 1 + \sinh^2 \pi] + 2 p_1 p_2 \cosh \pi \cosh([q_1 - q_2]_p - \pi) + 
p_2^2 [1 + \sinh^2 ([q_1 - q_2]_p - \pi)], \\
&q_2' = p_2^2 [ 1 + \sinh^2 \pi] + 2 p_1 p_2 \cosh \pi \cosh([q_2 - q_1]_p - \pi) + 
p_1^2 [1 + \sinh^2 ([q_2 - q_1]_p - \pi)], \\
&p_1' = - p_1 p_2 \sinh([q_1 - q_2]_p - \pi) [ p_1 \cosh \pi + p_2 \cosh([q_1-q_2]_p - \pi)], \\
&p_2' =- p_1 p_2 \sinh([q_2 - q_1]_p - \pi) [ p_2 \cosh \pi + p_1 \cosh([q_2-q_1]_p - \pi)]  .
\end{split}
\end{align}
Setting
\begin{align}
\notag
E(x) \doteq \frac{1}{\cosh \pi}  \cdot  \cosh([x]_p - \pi),
\quad
E'(x) \doteq \frac{1}{\cosh \pi}  \cdot \sinh([x]_p - \pi),
\end{align}
and using $q = q_2 - q_1$,
our system can be written in the more compact form 
\begin{align}
\label{per-pq-system-1}
\begin{split}
q_1' &=  \cosh^2 \pi \cdot ( p_1 + p_2E(q) )^2, \\
q_2' &=  \cosh^2 \pi \cdot  (  p_1 E(q) + p_2 )^2, \\ 
p_1' &=  \cosh^2 \pi \cdot p_1 p_2  ( p_1  + p_2 E(q) )  E'(q), \\
p_2' &=- \cosh^2 \pi  \cdot p_1 p_2   (p_1E(q) + p_2    ) E'(q)  . 
\end{split}
\end{align}
{\bf Initial Data.}
From this point, we make the same initial data assumptions as in the real-line case.
We take the positions the positions, $q_1$ and $q_2$  at time $t = 0$
to be
\begin{align}
\notag
q_1(0) = -a \quad \text{and} \quad q_2 = a, \quad \text{for some} \;\; a > 0.
\end{align}
For the initial momenta,  we shall assume  that at time $t=0$  that 
\begin{align}
\label{per-b-choice}
 p_2(0)=b\gg 1, \qquad    p_1(0)= -(b+\delta),\quad \delta>0.
\end{align}
With these assumptions, the initial profile $u_0(x)=u(x,0)$ is the  asymmetric periodic antipeakon-peakon
\begin{align}
\label{per-antpeak-peak-data}
u_0(x) 
&= 
-(b+\delta)  \cosh([x+a]_p - \pi)  + b \cosh([x-a]_p - \pi).
\end{align}
This initial profile for $u$ 
is displayed in the Figure 3.
\vskip0.2in
\hskip1in
\begin{minipage}{0.7\linewidth}
\hspace*{0cm}
\vspace*{0cm}
\begin{tikzpicture}[xscale=.8,yscale=.1]
%
%
\newcommand\X{7};
\newcommand\Y{2};
\newcommand\FX{11};
\newcommand\FY{11};
\newcommand\R{0.6};
\newcommand\B{2};
\newcommand\D{0.5};
\newcommand\A{1};
%
%
%
\draw[->,line width=1pt,black] (-4.5,0)--(7.5,0) 
node[above left] {\fontsize{\FX}{\FY}$x$};
\draw[->,line width=1pt,black] (0,-20)--(0,20) node[below right] {\fontsize{\FX}{\FY}$u_0$};
\draw[domain=(-4):(7),  variable=\x, 
red, line width=1.5pt, samples=500] 
plot ({\x},
        { -(\B + \D) * cosh((\x+\A) - (2*pi) * floor((\x+\A)/(2*pi)) - pi)      
	 +\B * cosh((\x-\A) - (2*pi) * floor((\x-\A)/(2*pi)) - pi)    
	    });
%
%
%
\draw[line width=1pt,black,dashed] 
({-(\A)},{-(cosh(pi)*(\B+\D) -4.5}) 
node[] { } 
node[below,xshift=-1.7cm] {\fontsize{\FX}{\FY}
$p_1(0)=-(b+\delta)\approx u_0(-a)$}
--
({-(\A)},0) 
node[above, xshift=-.21cm] {\fontsize{\FX}{\FY}$q_1(0)=-a$} 
node[] {$\bullet$} ;
\draw[line width=1pt,black,dashed] 
({(\A)},{cosh(pi)*(\B+\D)-11}) 
node[] {}
node[above,xshift=1.2cm] {\fontsize{\FX}{\FY}
$u_0(a)\approx b=p_2(0)$}
--
({(\A)},0) 
node[below, xshift=.2cm] {\fontsize{\FX}{\FY}$q_2(0)=a$} node[] {$\bullet$};
\draw[line width=1pt,black,dashed] 
({6.28},{0}) 
node[] { } 
node[above, xshift=-.21cm] {\fontsize{\FX}{\FY}
$2\pi$}
--
({6.28},0) 
node[] {$\bullet$} ;
\end{tikzpicture}
\end{minipage}
\vskip0.15in
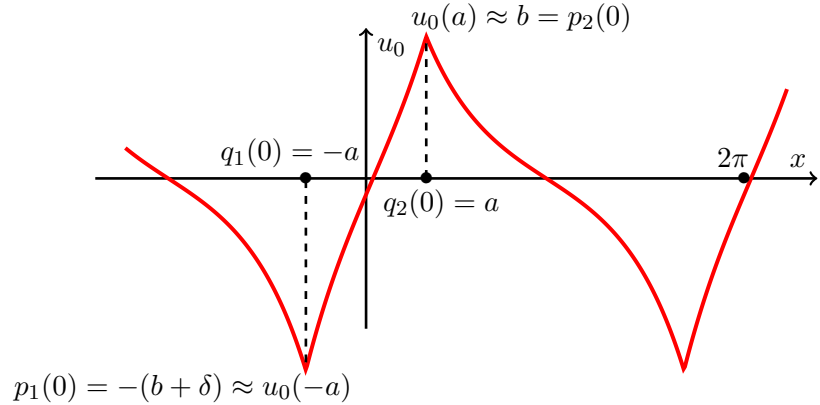
\centerline{Figure 3: Initial profile $u_0(x)$}
\vskip0.2in

Following the intuition we developed in the real-line case, we again will examine the 
ODE--system \eqref{per-pq-system-1} in the derived variables
$p,q,w,z$ given by
\begin{equation}
\label{per-q-p-w-def-periodic}
\begin{split}
q(t)&=q_2(t)-q_1(t), \quad \,\hspace{.02in} q(0)=2a>0,
\\
p(t)&=p_2(t)-p_1(t), \quad \, p(0)=2b+\delta>0,
\\
w(t)&=p_2(t)+p_1(t), \quad \hspace{-.01in} w(0)=-\delta<0,
\\
z(t)&=p_2(t)\cdot p_1(t), \quad \,\,\,\, z(0)=-b(b+\delta)<0.
\end{split}
\end{equation}

\vskip0.1in
\noindent
{\bf \large Deriving equations for $q$,
$p$, $w$ and $z$ on the circle.}
Beginning with $q$, we follow the same strategy as in the non-periodic case.  
We see that
\begin{align}
\notag
q'
&=
\cosh^2 \pi \cdot  (  p_1 E(q) + p_2 )^2
- 
\cosh^2 \pi \cdot  (  p_1 E(q) + p_2 )^2, \\
\notag
&=
\cosh^2 \pi \cdot 
(p_2-p_1) (p_2 +p_1)(1-E^2(q) ) , \\
\label{per-deriv-of-q=q2-q1}
&=
\cosh^2 \pi \cdot
p w(1-E^2(q)).
\end{align}
The computations for $p$, $w$ and $z$ follow the same strategy, 
and we arrive at the  $4 \times 4$ 
system 
\begin{equation}
\label{per-q-de-k=2}
\begin{split}
&q' \hspace{0.03in} = \cosh^2 \pi \cdot
p w(1-E^2(q)),  \quad \quad\quad \;\;  q(0)=2a>0,
\\  
&p' \hspace{0.03in} =   -  \cosh^2 \pi   \cdot w z (1+E(q)) E'(q) , \;\; \hspace{.025in}  p(0)=2b + \delta>0,
\\
&w' =  - \cosh^2 \pi \cdot  zp  (1- E(q)) E'(q)  , \; \hspace{.055in}    w(0)=-\delta<0,
\\
&z'  \hspace{0.03in} = \cosh^2 \pi  \cdot 
zwp E(q) E'(q),  \quad \quad\quad \;\; \hspace{.005in} z(0)=-b(b+ \delta)<0.
\end{split}
\end{equation}
This derived system of ODEs is more easily manipulated than the original 
$4\times4$ system, and we are now ready to tackle Proposition 1 in the periodic setting.

\vspace{.1in}


%
\begin{proposition}
[Periodic version of Proposition \ref{q-p-w-z-solutions}]
\label{per-q-p-w-z-solutions}
The system of differential equations  
\eqref{per-q-de-k=2}
has a unique smooth solution  $(q(t), p(t), w(t), z(t))$ in an  interval $[0, T)$,
for some $T>0$, such that  $z=z(t)$ is  decreasing
and in terms of $q$  is expressed  by the formula
\begin{align}
\label{per-z-q-form}
z
=
\frac{ - z_1}{\big(1 - E^2(q) \big)^{1/2}} < 0,
\quad
\text{where}
\quad
z_1=b(b+\delta) \big( 1 - E^2(q_0)  \big)^{1/2}>0,
\end{align}
$p=p(t)$ is  decreasing and as a function of $q$  is expressed  by the formula
\begin{align}
\label{per-p-q-form}
p(t)
=
\Big(
p_0^2 +
2
z_1 \Big[
\frac{1 + E(q(t)) }{\sqrt{ 1 - E^2(q(t)) }}
-
\frac{1 + E(q_0) }{\sqrt{ 1 - E^2(q_0)  }}
\Big]
\Big)^{1/2} > 0,
\end{align}
and
$w=w(t)$ is  decreasing and as a function of $q$  is expressed  by the formula
\begin{align}
\label{per-w-q-form}
w(t)
=
-\Big(
w_0^2 +
2z_1
\Big[
\frac{\sqrt{ 1 - E^2(q_0) }} { 1 + E(q_0) }
-
\frac{\sqrt{ 1 - E^2(q(t)) }} {1 + E(q(t)) }
\Big]
\Big)^{1/2} < 0.
\end{align}
The difference of the positions $q=q(t)$ is decreasing and satisfies the initial value problem
\begin{align}
\label{per-q-ivp}
&q'=-f(q)
\doteq
-
\cosh^2 \pi \cdot 
\Big(
w_0^2 +
2z_1
\Big[
\frac{\sqrt{ 1 - E^2(q_0) }} { 1 + E(q_0) }
-
\frac{\sqrt{ 1 - E^2(q(t)) }} {1 + E(q(t)) }
\Big]
\Big)^{1/2} 
\\&
\notag
\hskip1in
\cdot
\Big(
p_0^2 +
2
z_1 \Big[
\frac{1 + E(q(t)) }{\sqrt{ 1 - E^2(q(t)) }}
-
\frac{1 + E(q_0) }{\sqrt{ 1 - E^2(q_0)  }}
\Big]
\Big)^{1/2}
\cdot
(1 - E^2(q) ),
\\&
q(0)=q_0=2a>0.
\notag
\end{align}
Furthermore,  the initial value problem  \eqref{per-q-ivp}
  is dominated by the simpler initial value problem 
\begin{align}
\notag
q'=-g(q)
\doteq
-q_1 \big(1-e^{-2q}\big)^{3/4},
\quad
0<q(0)=2a<1/2,
\end{align}
where 
\begin{equation}
\notag
q_1
=
\delta 
\sqrt{2b(b+\delta)}\cdot q_0^{1/4}.
 \end{equation} 
\end{proposition}
%
%


{\bf Proof.}  
We begin by solving for $p,w$ and $z$ in terms of $q$. After this task is completed, 
we can form an autonomous equation for $q$ by substituting in these results.


\vspace{.1in}

{\bf Expressing $z$ in terms of $q$.} Using the equation for $z'$ and $q'$ we find 
\begin{align}
\notag
&\frac{z'}{q'}=\frac{ \cosh^2 \pi  \cdot  zwp E(q) E'(q)  }{  \cosh^2 \pi \cdot p w(1-E^2(q))  }
=
\frac{z E(q) E'(q)}{1 - E^2(q)} \quad \text{or} \quad
\frac{z'}{z}=\frac{ E(q) E'(q) q'}{1 -  E^2(q) }.
\end{align}
Since $z(0)<0$ we shall  assume that $z(t)$ will remain negative.  
Therefore, from the last relation we have 
\begin{align}
\notag
\frac{d}{dt}[\ln (-z)]
=
-\frac 12 \frac{d}{dt}[\ln ( 1 - E^2 (q)  )].
\end{align}
Integrating this equation from $0$ to $t$ gives
\begin{align}
\notag
\ln\Big[\frac{z(t)}{z_0}\Big]
=
-\frac 12 \ln \Big[\frac{1 - E^2(q(t)) }{ 1  -  E^2(q_0) }\Big].
\end{align}
Solving for $z(t)$, we find  formula \eqref{per-z-q-form}
for $z$ in terms of $q$.


\vspace{.1in}

{\bf Expressing $w$ in terms of $q$.} Dividing the 
equation for $w'$ by the equation for $q'$ we have 
\begin{align}
&\frac{w'}{q'}
=
\frac{  - \cosh^2 \pi \cdot  zp  (1- E(q)) E'(q) }{ \cosh^2 \pi \cdot p w(1-E^2(q)) }
=
\frac{  - z  (1- E(q)) E'(q) }{ w(1-E^2(q)) }, \\
&
\text{or}
\quad
ww'=- z\cdot \frac{   ( 1 - E(q) ) E'(q) q'}{  1 - E^2(q)  }.
\end{align}
Substituting the formula for $z$ given by \eqref{per-z-q-form} into the above equation
gives us
\begin{align}
\label{per-ww'-eq}
ww'
=
\frac{z_1}{\big(1 - E^2(q) \big)^{1/2}} 
\cdot \frac{   (1 - E(q) ) E'(q) q'}{1 -  E^2(q)  }
=
\frac{z_1(1 -  E( q) )  E'(q) q'}{\big(1 - E^2(q)  \big)^{3/2}}.
\end{align}
Making the change of variables
$u=E(q(t))$, $du=E'(q(t)) q'(t) dt$, we obtain
\begin{align*}
&\int
\frac{(1  - E(q) ) E'(q) q'}{\big( 1  -  E^2(q) \big)^{3/2}}
dt
=
\int \frac{1 - u }{\big(1 - u^2 \big)^{3/2}} du
=
-
\frac{ \sqrt{1 - u^2}}{ 1 + u } + C
=
-
\frac{ \sqrt{ 1 - E^2(q(t))}}{1 + E(q(t))} + C.
\end{align*}
Therefore, relation \eqref{per-ww'-eq} reads as 
\begin{align}
\frac{d}{dt}\Big[\frac12 w^2\Big]
=
-z_1 \frac{d}{dt}\Big[
\frac{\sqrt{ 1 - E^2(q(t)) }} { 1+  E(q(t)) }
\Big].
\end{align}
Integrating this equation from $0$ to $t$ gives
\begin{align}
\frac12 \Big[w^2(t)-w_0^2\Big]
=
z_1 \Big[
\frac{\sqrt{ 1 - E^2(q_0) }} { 1 + E(q_0) }
-
\frac{\sqrt{ 1 - E^2(q(t)) }} { 1 + E(q(t)) }
\Big].
\end{align}
We are thus able to solve for $w(t)$ in terms of $q(t)$,
which gives us  formula \eqref{per-w-q-form}.

\vspace{.1in}

{\bf Expressing $p$ in terms of $q$.} Dividing the 
equation for $p'$ by the equation for $q'$ we have 
\begin{align}
&\frac{p'}{q'}
=
\frac{  -  \cosh^2 \pi   \cdot w z (1+E(q)) E'(q)) }{  \cosh^2 \pi \cdot p w(1-E^2(q))}
=
\frac{  -     z (1+E(q)) E'(q)) }{   p (1-E^2(q))}, \\
&
\text{or}
\quad
pp'= - z\cdot \frac{ (1 + E(q) ) E'(q)  q'}{ 1 -  E^2(q) }.
\end{align}
Substituting in the above relation the formula for $z$
given by \eqref{per-z-q-form} we have
\begin{align}
\label{per-pp'-eq}
pp'
=
\frac{z_1}{\big( 1 - E^2(q)  \big)^{1/2}} \cdot \frac{  ( 1 + E(q) ) E'(q) q'}{1 - E^2(q) }
=
\frac{z_1( 1 + E(q) ) E'(q) q'}{\big( 1 - E^2 (q)   \big)^{3/2}}.
\end{align}
Next, we make  the change of variables
$u=E(q(t))$,  $du=E'(q(t)) q'(t) dt$  and get 
\begin{align*}
\int
\frac{(1 + E(q) ) E'(q) q'}{\big( 1 -  E^2(q)  \big)^{3/2}}
dt
=
\int \frac{1 + u }{\big( 1 -  u^2 \big)^{3/2}} du
=
\frac{1 + u }{\big( 1 -  u^2 \big)^{1/2}}+C
=
\frac{1 + E(q) }{\big( 1 - E^2(q) \big)^{1/2}} +C.
\end{align*}
Therefore, relation \eqref{per-pp'-eq} reads as follows
\begin{align}
\notag
\frac{d}{dt}\Big[\frac12 p^2\Big]
=
z_1 \frac{d}{dt}\Big[\frac{1 + E(q)  }{\sqrt{1 - E^2(q)  }}\Big].
\end{align}
Integrating this equation from $0$ to $t$ gives us
\begin{align}
\notag
\frac12 \Big[p^2(t)-p_0^2\Big]
=
z_1 \Big[
\frac{1 + E(q(t)) }{\sqrt{1 - E^2(q(t)) }}
-
\frac{1 + E(q_0) }{\sqrt{ 1 - E^2(q_0)  }}
\Big],
\end{align}
and we are able to solve for $p(t)$ 
and obtain formula \eqref{per-p-q-form}.

%
%
%
%

\vskip0.1in

{\bf Solving the $q$ ODE.}
Starting with the differential equation for $q$, which is
$q' = \cosh^2 \pi \cdot
p w(1-E^2(q))$, we substitute in for $w$ and
$p$ their expressions \eqref{per-w-q-form} and \eqref{per-p-q-form}  respectively.
We consequently obtain the following autonomous differential 
equation for $q$
\begin{align}
\label{per-q-autonom-de}
&q'=-f(q)
\doteq
\cosh^2 \pi \cdot \Big\{
-\Big(
w_0^2 +
2z_1
\Big[
\frac{\sqrt{ 1 - E^2(q_0) }} { 1 + E(q_0) }
-
\frac{\sqrt{ 1 - E^2(q(t)) }} {1 + E(q(t)) }
\Big]
\Big)^{1/2} \Big\}
\\&
\notag
\hskip1in
\cdot
\Big(
p_0^2 +
2
z_1 \Big[
\frac{1 + E(q(t)) }{\sqrt{ 1 - E^2(q(t)) }}
-
\frac{1 + E(q_0) }{\sqrt{ 1 - E^2(q_0)  }}
\Big]
\Big)^{1/2}
\cdot
(1 - E^2(q) ),
\\&
q(0)=q_0=2a>0.
\notag
\end{align}
Next, we observe that  
\begin{align}
\label{per-w-simplification}
\frac{\sqrt{ 1 - E^2(q_0) }} {1 +  E(q_0) }
-
\frac{\sqrt{ 1 - E^2(q(t)) }} {1 +  E(q(t)) }
\ge 0, 
\quad
0\le q\le q_0 < \pi.
\end{align}
This inequality follows from the fact that 
\begin{align}
\notag
\Big(
\frac{\sqrt{1 - E^2(x)}}{1 + E(x)}
\Big)'
=
\frac{-E'(x)}{2\sqrt{1+E(x)}}.
\end{align}
Since the denominator is always positive, the sign of this derivative is controlled by numerator,
$-E'(x) = -\frac{1}{\cosh (\pi)} \sinh ([x]_p - \pi)$,
which is positive for $x \in [0,\pi)$.
Next, we have
\begin{align}
\label{per-p-simplification}
p_0^2  -
2z_1 \frac{1 + E(q_0) }{\sqrt{ 1 - E^2(q_0) }}
\ge 0
&\iff
\frac{(2b + \delta)^2}{ 2b(b+\delta)}\ge 1 +  E(q_0).
\end{align}
Our choice of initial data allows for the inequality $1 + E(q_0) \le 2$, and 
we have
\begin{align}
\notag
\frac{(2b+\delta)^2}{ 2b(b+\delta)}\ge 2
\iff
4b^2+4b\delta+\delta^2>4b^2+4b\delta
\iff
\delta^2>0,
\quad
\text{which is true.}
\end{align}

\vskip0.1in
Now, using  \eqref{per-w-simplification} and  \eqref{per-p-simplification} 
we see that  the function $f(q)$ in the right-hand side  of the differential equation \eqref{per-q-autonom-de} can be bounded from below 
as follows 

\begin{align}
\notag
f(q)
&\ge
\cosh^2 \pi \cdot
\delta
\Big(
2 b(b+\delta) \big(1-E^2 (q_0) \big)^{1/2}\Big)^{\frac 12}
\Big(
 \Big[
\frac{1 + E(q(t)) }{\sqrt{1  - E^2(q(t)) }}
\Big]
\Big)^{\frac 12}
\cdot
(1 - E^2(q) ).
\end{align}
To continue in our objective of finding a simpler dominating function for $f$, 
analogous to the strategy in the real-line case of this proof, 
we use the fact that $E(q) \ge 0$  in conjunction with the following 
lemma.
\begin{lemma}
\label{per-E-inequality}
For $c \ge 2\cosh^2(\pi)/\sinh(2 \pi - 1)$,  and $x \in [0,1/2]$,  
\begin{align}
\label{per-E-inequality-1}
c(1 - E^2(x)) \ge 1 - e^{-2x}.
\end{align}
Furthermore, we have the inequality
\begin{align}
\label{per-E-inequality-2}
1 - E^2(q_0) \ge \frac{1}{3} q_0.
\end{align}
\end{lemma}
In particular, we will take $c = 3$ in later computations.

\vspace{.1in}
{\bf Proof.}  Define the function
\begin{align}
\notag
f(x) \doteq \Big( 1 - e^{-2x} \Big) - \Big( c \cdot [1 - E^2(x) ] \Big),
\end{align}
Computing the derivative of $f(x)$ shows that it will be negative for
$x \in (0,1/2]$,
and
\begin{align}
\notag
 c  \ge  \frac{2 \cosh^2 \pi}{\sinh(2 \pi - 1)}.
\end{align}
As the \eqref{per-E-inequality-1} has been established, we now
move onto proving \eqref{per-E-inequality-2}.
This inequality is obtained by  applying  our first inequality
and then using the exponential 
inequality.
We get
\begin{align}
\notag
1 - E^2(q_0) \ge \frac{1}{3} (1 - e^{-2q_0}) \ge \frac{1}{3} q_0.
\qed
\end{align}

\vspace{.1in}

With the above lemma, we are now ready to return to the proof of the proposition.  
\vspace{.1in}

{\bf Dominating Equation (periodic version).}
Using the above inequalities, and following the same strategy as in the non-periodic case, we obtain
\begin{align}
\notag
f(q)
&\ge
\Big[ \frac{  \sqrt{2}  \cosh^2 \pi}{3} \cdot
\delta
\cdot
\sqrt{ b(b+\delta) } \cdot q_0^{1/4}  
 \Big] 
\cdot \Big(1 - e^{-2q} \Big)^{3/4}.
\end{align}
Since 
$\frac{  \sqrt{2}  \cosh^2 \pi}{3} \ge 1,$
we can remove this factor as we are bounding from below.   Consequently, $f(q)$ has precisely the same lower bound as in the real-case give by
\begin{align}
\notag
f(q)
&\ge
\delta
\cdot
\sqrt{ b(b+\delta) } \cdot q_0^{1/4}  
\cdot \Big(1 - e^{-2q} \Big)^{3/4}
=q_1  \Big(1 - e^{-2q} \Big)^{3/4},
\end{align}
where the constant $q_1$ is given by
\begin{align}
\notag
q_1
=
\delta
\cdot
\sqrt{  b (b+\delta)} \cdot q_0^{1/4}.
\end{align}
Thus, we see that the complicated initial value problem  for $q$
\eqref{per-q-autonom-de} is dominated by 
\begin{align}
\label{per-q-autonom-de-simpler}
q'
=
- q_1 \big(1- e^{-2q}  \big)^{3/4},
\quad
q(0)=q_0=2a>0.
\end{align}
This ODE is precisely the same as the one derived in the real-line case.  
Therefore we can immediately arrive at the same conclusions for $q$.  
\begin{proposition}[Periodic version of Proposition \ref{zero-of-q-lifespan}]
\label{per-zero-of-q-lifespan}
If $r<1$ then for given $q_0\in (0, 1/2)$ and $q_1>0$
the solutions to the initial value problem
\begin{align}
\label{per-q-ivp-r}
\frac{dq}{dt}  
=
-g_r(q)
\doteq
 -q_1
 \left(1-e^{-2q}\right)^r,
 \quad
 q(0)=q_0,
 \end{align}
which begins positive and is decreasing, becomes zero in finite time
$T$ given by 
\begin{align}
\label{per-T-lifespan-r}
T
=
\int_0^{q_0}
 \frac{dq}{g_r(q)}
 =
\frac{1}{q_1}
\int_0^{q_0}
 \frac{dq}{\left(1-e^{-2q}\right)^r}
\simeq
 \frac{1}{1-r}\, \frac{q_0^{1-r}}{q_1}.
\end{align}
\end{proposition}

\begin{corollary} [Periodic version of Corollary \ref{zero-of-q-lifespan-of-u}]
\label{per-zero-of-q-lifespan-of-u}
If $0<q_0<1/2$ and $b>1$,  $\delta >0$ satisfy condition 
\eqref{per-b-choice} 
then the solution  to the initial value problem
\eqref {per-q-ivp}
begins positive, is decreasing, and becomes zero in finite time
$T$ given  by
{\small
\begin{align}
\label{per-T-lifespan-3/4}
T
=
\int_0^{q_0}
 \frac{dq}{f(q)}
 \le
 \frac{1}{q_1}
\int_0^{q_0}
 \frac{dq}{\left(1-e^{-2q}\right)^{3/4}}
\simeq
 \frac{q_0^{1/4}}{q_1}
 \simeq
  \frac{q_0^{1/4}}{
\delta
\sqrt{2b(b+\delta)}\cdot q_0^{1/4}
  }
   \simeq
  \frac{1}{
\delta
\sqrt{2b(b+\delta)}
  }.
  \end{align}
  }
\end{corollary}

\vspace{.1in}
We summarize the above results  
in the following  Theorem.
\begin{theorem}[Periodic version of Theorem
 \ref{cor:NE-2-peakon}]
\label{per-cor:NE-2-peakon}
For given $0<a\le 1/4$, $b > 1$  and $\delta > 0 $ satisfying  condition \eqref{per-b-choice},  
the initial value problem for the positions $q_1, q_2$ and the momenta $p_1, p_2$
\begin{align}
\label{per-pq-system-2}
\begin{split}
q_1' &=  \cosh^2 \pi \cdot ( p_1 + p_2E(q) )^2,
\quad\quad\quad\quad \;\;\;\;\; \; \hspace{.015in}  q_1(0) = -a,
\\
q_2' &=  \cosh^2 \pi \cdot  (  p_1 E(q) + p_2 )^2, 
\quad\quad\quad\quad \;\;\;\;\;\; \hspace{.015in} q_2(0) = a,
\\ 
p_1' &=  \cosh^2 \pi \cdot p_1 p_2  (p_1  + p_2 E(q) )   E'(q), 
\quad\quad p_1(0) = -(b + \delta),
\\
p_2' &=- \cosh^2 \pi  \cdot p_1 p_2  ( p_1E(q) + p_2  )  E'(q),
\quad p_2(0) = b,
\end{split}
\end{align}
has a unique smooth solution $(q_1(t), q_2(t), p_1(t), p_2(t))$
with a finite lifespan  $T$, which is the zero of $q=q_2-q_1$, satisfying the estimate 
\eqref{per-T-lifespan-3/4}
and such that 
\begin{align}
\notag
&p_1=\frac{w-p}{2}<0, \,\, decreasing, \,\, 
 \lim_{t \to T^{-}} p_1(t) = -\infty, \,\, 
  \text{and}  \,\,  -p_1\simeq p\simeq q^{-1/4}, \;\; \text{and} \\ 
\notag
&p_2=\frac{w+p}{2}>0, \,\, increasing, \,\, 
 \lim_{t \to T^{-}} p_2(t) = \infty, \,\, 
  \text{and}  \,\,  p_2\simeq p  \simeq q^{-1/4},
\end{align}
where $p$ and  $w$ are given in Proposition \ref{per-q-p-w-z-solutions}.
Also,  $w=p_1+p_2$ is decreasing from 
$ w_0 < 0$ to $w_T$, where
\begin{align}
\notag
&w_T
\doteq
 \lim_{t \to T^{-}} w(t)
 =
-\Big(
\delta^2 +
2   b(b+\delta)
(1-E(2a))
\Big)^{\frac 12}.
\end{align}
Finally,  the  2-peakon  
\begin{align}
\notag
u(x,t) = p_1(t) \cosh([x-q_1(t)]_p - \pi)  + p_2(t) \cosh([x-q_2(t)]_p - \pi),
\, x\in \tor, \, 0<t<T.
\end{align}
is a solution to NE with following the asymmetric  antipeakon-peakon initial profile
\begin{align}
\notag
u(x,0) 
= 
-(b+\delta) 
\cosh([x+a]_p - \pi)
+b \cosh([x-a]_p - \pi),
\end{align}
The quantities $p_1$, $p_2$, $q$ and $w$ have similar properties
to their analogues defined on the line, and we refer to Figure 2 for a visualization of them. 
\end{theorem}
%

%

%
%
%
%
%
%
\subsection{Calculating the Norm on the circle}
\label{subsec:3}
We begin with the following proposition which summarizes the calculation of the $H^s$ norm of $u$.  This computation is nearly identical to non-periodic case with the exception of an extra
factor of $\sinh^2 \pi$.

\begin{proposition}[Periodic version of Proposition \ref{norm-estimate-NE}]
\label{per-norm-estimate-NE}
Let $u(t)$ be the two-peakon solution \eqref{per-NE-2-peakon} to the NE equation. Then on $[0,T)$ we have 
\begin{align}
\notag
&\|u(t)\|_{H^s}^2
=
16 \sinh^2 \pi \cdot r(t) p_1^2(t)     Q_s(q)
+
4 \sinh^2 \pi \cdot c_s\big(1-r(t)\big)^2 p_1^2(t), \\
\notag
&
\text{with}
\quad 
r(t)\doteq -\frac{p_2(t)}{p_1(t)},
\end{align}
where    $c_s= \sum_{-\infty}^\infty  (1+n^2)^{s-2}  $  and  
$Q_s(q)$, which  is given below,
satisfies  the estimates:
\begin{align}
\notag
 Q_s(q)
 \doteq
\sum_{n=-\infty}^\infty
(1+n^2)^{s-2} \sin^2 \bigg(\frac{q n}{2} \bigg) 
\simeq
\begin{cases}
q^{3-2s},     \quad \;\quad\quad 1/2 < s < 3/2, \\
q^2 \cdot \ln(1/q),   \;\;\;\, s = 1/2, \\
q^2,\quad \quad \quad \quad \quad s < 1/2.
\end{cases}
\end{align}
\end{proposition}

{\bf Proof.} 
We begin by noting that the Fourier transform of $E$ is calculated as
\begin{align}
\notag
\widehat{E}(n)
&= 
\Big(2  \cdot  \frac{\sinh(\pi)}{\cosh(\pi)} \Big) \cdot \frac{1}{1 + n^2}. 
\end{align} 
Recalling that the 2-peakon $u$ can be written as
\begin{align}
\notag
u(x,t) 
&=
\cosh \pi \cdot \Big(
p_1(t) E(x - q_1(t)) + p_2(t) E(x-q_2(t)) 
\Big),
\end{align}
we can express the Fourier transform of $u$ as
\begin{align}
\notag
\widehat{u}(n,t)
 &=
\frac{2  \sinh \pi }{1+n^2}  \cdot 
 p_1 e^{-i n  q_1} \cdot \Big(1 + \frac{p_2}{p_1} e^{-i n  q} \Big).
 \notag
\end{align}
Taking the square of the $H^s$ norm of this quantity, we obtain
\begin{align}
\label{per-norm-of-u-1}
\|u(t)\|_{H^s}^2
&=
\sum_{n = -\infty}^\infty
(1 + n^2)^s |\widehat{u}(n,t)|^2
=
4 \sinh^2 \pi  \cdot p_1^2 
\sum_{n = -\infty}^\infty
(1+n^2)^{s-2} \Big| 1 + \frac{p_2}{p_1} e^{-i n q} \Big|^2.
\end{align}
Using Proposition \ref{per-q-p-w-z-solutions} we see that
\begin{align}
\label{per-q1}
r=r(t)\doteq -\frac{p_2(t)}{p_1(t)}=\frac{p+w}{p-w}<1,
\quad
\text{and}
\quad
 r(t) \nearrow 1 \text{ as } t \nearrow T.
\end{align}
Using  $r$ we write \eqref{per-norm-of-u-1} as follows
\begin{align}
\label{per-u-norm-1}
\|u(t)\|_{H^s}^2
=
4 \sinh^2 \pi \cdot  p_1^2 
\sum_{n = -\infty}^\infty
(1+n^2)^{s-2} \Big| 1-r e^{-i n q}  \Big|^2.
\end{align}
Expanding out the square of the absolute value inside of the sum 
\eqref{per-u-norm-1}, we have
\begin{align}
\notag
|re^{i q n} - 1|^2
&=
(1-r)^2+ 4 r \sin^2 \bigg(\frac{qn}{2} \bigg).
\end{align}
We therefore obtain the formula 
\begin{align}
\notag
\|u(t)\|_{H^s}^2
=
16 \sinh^2 \pi \cdot  r p_1^2     Q_s(q)
+
4   \sinh^2 \pi \cdot  c_s(1-r)^2 p_1^2,
\end{align}
where 
\begin{align}
\notag
c_s=  \sum_{n = -\infty}^\infty (1+n^2)^{s-2}
\quad \text{and} \quad
Q_s(q)
=
\sum_{n = -\infty}^\infty
(1+n^2)^{s-2} \sin^2 \bigg(\frac{qn}{2} \bigg) .
\end{align}
From this point, we note that $Q_s$ has already been
estimated in this periodic setting in \cite{hhg}.  Using 4.25 from \cite{hhg}, and noting 4.28, where the norm is expanded into the sum of the squares of sines, we have  
\begin{align}
\label{per-Qs}
Q_s 
\simeq
\sum_{n = 1}^\infty \sin^2 \Big(\frac{q n}{2} \Big) (1 + n^2)^{s-2}
\simeq
\begin{cases}
q^{3/2 - s}, & 1/2 < s < 3/2, \\
q \sqrt{\ln(1/q)}, & s = 1,2, \\
q, & s < 1/2.   \qquad \qed
\end{cases}
\end{align}
%
%

%
%
%
%
%
%
%

\subsection{Small lifespan and initial data on the circle}
\label{subsec:4}
This section follows the same argument as in the real-line case, with the exception of an 
extra factor of $\sinh^2 \pi$ stemming from the periodic version of the norm-estimates.  
We begin by assuming that 
\begin{align}
\notag
 p_2(0)=b\gg1 \text{ and }-p_1(0)=  b+\delta,\, \delta>0,
\end{align}
so that  the conditions for the existence
of our 2-peakon with the lifespan estimate
\eqref{per-T-lifespan-3/4}
hold. 

\medno
{\bf Lifespan Estimate.}   For given $\ee>0$,
we need to find $b>1$ such that $T<\ee$. 
Since, by Proposition \ref{per-q-p-w-z-solutions} we have
\begin{align}
\notag
T
\lesssim
\frac{1}{
\delta 
\sqrt{2b(b+\delta)}}
\le
\frac{1}{\delta b},
\end{align}
we must have
\begin{align}
\label{per-T-small-1}
\frac{1}{
\delta b}
\le \ee \iff  b\ge \delta^{-1}  \ee^{-1}.
\end{align}

\medno
{\bf Initial Data Estimate.} Now, for the same $\ee>0$
we need to find  $q_0<1/8$ such that
$\|u_0\|_{H^s} <\ee$. 
For this we use Proposition \ref{per-norm-estimate-NE},
from which we have, recalling that $r(t)= -\frac{p_2(t)}{p_1(t)}$, 
\begin{align}
\notag
\|u(0)\|_{H^s}^2
=
16 \sinh^2 \pi \cdot b(b+\delta)    Q_s(q_0)
+
4 \sinh^2 \pi \cdot c_s\delta^2.
\end{align}
This identity implies 
\begin{align}
\notag
\|u(0)\|_{H^s}^2
\le
32 \sinh^2 \pi \cdot b^2    Q_s(q_0)  
+
4 \sinh^2 \pi \cdot c_s\delta^2.
\end{align}
%
%
%
%
%
%
{\bf Case $1/2 < s < 3/2$:} Then  by Proposition \ref{per-norm-estimate-NE}
we have  $Q_s(q_0) \lesssim q_0^{3-2s}$ and therefore
\begin{align}
\notag
\|u(0)\|_{H^s}^2
\le
C_s  \sinh^2 \pi \cdot b^2   q_0^{3-2s}
+
4  \sinh^2 \pi \cdot c_s\delta^2.
\end{align}
For having  $\|u_0\|_{H^s} <\ee$ it suffices 
to choose $q_0$ and $\delta$ such that
$
C_s \sinh^2 \pi \cdot b^2   q_0^{3-2s}
+
4 \sinh^2 \pi \cdot c_s\delta^2
\le \ee^2$, or
\begin{align}
\notag
4\sinh^2 \pi \cdot c_s\delta^2  \le \frac{\ee^2}{2}
\quad
\text{ and }
\quad
C_s \sinh^2 \pi \cdot b^2   q_0^{3-2s} \le \frac{\ee^2}{2}.
\end{align}
The first  inequality holds if   
\begin{align}
\label{per-delta-choice}
\delta \le \frac{\ee}{2 \sinh \pi \sqrt{2c_s }}.
\end{align}
Taking into consideration \eqref{per-delta-choice} and 
\eqref{per-T-small-1}, the second  inequality holds if 
\begin{align}
\notag
q_0^{3-2s} \le \frac{\ee^2}{2C_s  \sinh^2 \pi \cdot b^2 }
\le
\frac{\ee^2}{2C_s \delta^{-2} \ee^{-2} }
=
\frac{\delta^{2}\ee^4}{2C_s}
\le
\frac{\ee^{2}\ee^4}{8c_s\cdot 2C_s},
\end{align}
or 
\begin{align}
\notag
q_0 
\le
\Big(
\frac{\ee^6}{16c_sC_s}
\Big)^\frac{1}{3-2s}.
\end{align}
%
%
%
%
{\bf Case $s \le 1/2$:} For a such Sobolev exponent $s$
we have   $\|u(0)\|_{H^s} \le \|u(0)\|_{H^1}$. 
This combined with  Proposition \ref{per-norm-estimate-NE},
which tells us that   $Q_1(q_0) \lesssim q_0$,
gives
\begin{align}
\notag
\|u(0)\|_{H^s}^2
\le
\|u(0)\|_{H^1}^2
\le
C_1b^2   q_0
+
4c_1\delta^2.
\end{align}
Thus  $\|u_0\|_{H^s} <\ee$  if 
$q_0$ and $\delta$ satisfy the inequalities 
$
4c_1\delta^2  \le \frac{\ee^2}{2}
$
 and 
$
C_1b^2   q_0 \le \frac{\ee^2}{2}.
$
These inequality holds if 
\begin{align}
\notag
\delta \le \frac{\ee}{2\sqrt{2c_1}}
\quad
\text{and}
\quad
 q_0 
 \le
\frac{\ee^6}{16c_sC_s}.     
\end{align}
%

%
%

\subsection{Norm-Inflation and illposedness on the circle}
\label{subsec:5}
%
%
%

%
%
From Proposition \ref{per-norm-estimate-NE}
we have
\begin{align}
\label{per-u-norm-NE-fin}
\|u(t)\|_{H^s}^2
=
16 \sinh^2 \pi \cdot  r(t) p_1^2(t)     Q_s(q)
+
4 \sinh^2 \pi \cdot   c_sp_1^2(t)\big(1-r(t)\big)^2,
\end{align}
We see that the same argument as in Section \ref{sec:5} holds, with the simple inclusion
of a factor of $\sinh^2\pi$.  Following these arguments, we see that
\begin{equation}
\label{per-limit-H^s-norm-u}
\lim_{t\to T} \|u(t)\|_{H^s}^2
=
\begin{cases}
\infty\quad
\text{(inflation)},     &    5/4 < s < 3/2, \\
\text{may not exist},  &  s = 5/4, \\
4 \sinh^2 \pi \cdot c_s\Big[
\delta^2 +
2b(b+\delta)
\cdot
(1-e^{-q_0})\Big],  &  s < 5/4.
\end{cases}
\end{equation}
Therefore when   $5/4<s <3/2$ we have norm inflation
and ill-posedness for the Novikov equation.  \qed

%
%
%
%
%
%
%
\subsection{Non-Uniqueness for $s<5/4$ on the circle}
\label{subsec:6}

As in the case on the line, the NE admits non-unique solutions once we take the Sobolev
exponent $s < 5/4$.  This is an equally interesting result as the periodic 2-peakons 
maintain the same collision properties as non-periodic ones. 

\begin{theorem}[Nonuniqueness - Periodic version of Theorem \ref{NE-non-uniqueness}] 
\label{per-NE-non-uniqueness}
For $s < 5/4$  NE  admits non-unique solutions.
\end{theorem}

Our proof of non-uniqueness in the periodic setting
 again follows the same strategy used in the real-line case.
We again examine the behavior of the limit as 
$t \to T^-$ of the 2-peakon
solution $u$ with initial data given in \eqref{per-antpeak-peak-data}.  Once this limit
has been established in the desired ways, the same argument as in the real-line 
case implies non-uniqueness.

\begin{proposition}  
[Pointwise limit - Periodic version of Proposition \ref{u-pointwise-limit}]
\label{per-pointwise-limit} 
For each $x \in \rr$ we have 
\begin{align}
\label{per-sol-t-limit}
\lim_{t \to T} u(x,t) = w_T \cosh([x - q_T]_p - \pi)  \doteq v_T(x). 
\end{align}
where
\begin{align}
& q_T \doteq \lim_{t \to T^-} q_1(t) = \lim_{t \to T^-}  q_2(t) 
\quad \text{and} \quad
w_T \doteq \lim_{t \to T^-} w(t).
\end{align}
\end{proposition}

\vspace{.1in} 
{\bf Remark.} We can avoid the multiple cases needed in the real-line version of this proof as we do not need to expand out an absolute value.  Here, as we are using the hyperbolic cosine, we will have both $e^{x}$ and $e^{-x}$ present thus avoiding the need to break into cases. 

\vspace{.1in}
{\bf Proof.}
Our solution $u$ is a  $2\pi$-periodic function, and we will restrict our attention to the interval $[0,2\pi]$.  As we know that the limits of $q_1$ and $q_2$ exists, we will further restrict our attention to after some time $t_0 >0 $ such that these position function remain within a single period.  This will avoid any complications of moving between periods which will require using the floor function in our definition.    Using the exponential definition of the hyperbolic cosine, we get
\begin{align}
\notag
u(x,t)
&=
\frac{1}{2}
\Big[ 
e^{\pi - x}( p_1 e^{q_1}+ p_2 e^{q_2}) + e^{-\pi +x}(p_1 e^{-q_1} + p_2 e^{-q_2} )
\Big].
\end{align}
Rewriting this expression to generate terms containing $w$ gives us
\begin{align}
\label{per-ptwise-per-limit-2}
u(x,t)
&=
\frac{1}{2}
\Big[
e^{\pi - x}
\Big(e^{q_2} p_1  (e^{-q} - 1)  + w e^{q_2} \Big)
+
e^{-\pi + x}
\Big(
w e^{-q_1} + e^{-q_1} p_2 (e^{-q} - 1)
\Big)
\Big].
\end{align}
Taking the limit as $t \to T^-$ of \eqref{per-ptwise-per-limit-2}, and using the limit 
established in Lemma \ref{lem-p-pointwise-limits}, 
we obtain
\begin{align}
\notag
\lim_{t \to T^-}
u(x,t)
&=
w_T \cosh([x- q_T]_p - \pi).
\qed
\end{align}

We next
demonstrate that $u$ converges to $v_T$ as $t \to T^-$
in $L^r$.  
\begin{proposition}[Convergence in $L^r$  - Periodic version of Proposition \ref{Lr-conv}] 
\label{per-Lr-conv}
For our antipeakon-peakon solution $u$ to NE, we have  
\begin{align}
\lim_{t \to T^-} \| u(x,t) - v_T(x) \|_{L^r} = 0.
\end{align}
\end{proposition}

{\bf Proof.}   The same remarks that we made for the pointwise proof apply here as 
to taking a $t_0 > 0$ such that $q_1$ and $q_2$ lie within a single $2\pi$ period after
time $t_0$.
Analogous to the pointwise limit, as we have both $e^{x}$ and $e^{-x}$ present 
in our hyperbolic cosines, we will not have to break our argument into cases in order to simplify the
absolute values.  This fact also allows us to bypass the restriction $1 \le r <4$ as we do not
cut the domain of the integration, creating the situation we saw in the real-line case on the sub-integral
on $[q_1,q_2]$.  After rewriting the hyperbolic cosines in their exponential form we get
\begin{align*}
\| u(x,t) - v_T (x) \|_{L^r}^r
&=
\frac{1}{2}
\int_0^{2\pi}
\Big|
e^{x - \pi}
\Big( p_1 e^{-q_1} + p_2 e^{-q_2} - w_T e^{-q_T} \Big)   \\
&\quad\quad\quad\quad\quad+
e^{\pi - x} \Big(
p_1 e^{q_1} + p_2 e^{q_2} - w_T e^{q_T} \Big) \Big|^r dx.
\end{align*}
Using Jensen's inequality, and evaluating the resulting integrals, we get
\begin{align}
\notag
&
\frac{1}{2} \int_0^{2\pi}
\Big|
e^{x - \pi}
\Big( p_1 e^{-q_1} + p_2 e^{-q_2} - w_T e^{-q_T} \Big) 
+
e^{\pi - x} \Big(
p_1 e^{q_1} + p_2 e^{q_2} - w_T e^{q_T} \Big) \Big|^r dx \\ \notag
&\quad
\le
\frac{e^{r\pi} - e^{-r \pi}}{r}   \cdot  \Big| p_1 e^{-q_1} + p_2 e^{-q_2} - w_T e^{-q_T} \Big|^r
+
\frac{e^{r\pi} - e^{- r \pi}}{r} \cdot
\Big| p_1 e^{q_1} + p_2 e^{q_2} - w_T e^{q_T} \Big|^r.
\end{align}
Using  Lemma \ref{lem-p-pointwise-limits},
we have
\begin{align}
\notag
&\lim_{t \to T^-}\Big| p_1 e^{q_1} + p_2 e^{q_2} - w_T e^{q_T} \Big|^r = 0
\quad \text{and} \quad
\lim_{t \to T^-}\Big| p_1 e^{q_1} + p_2 e^{q_2} - w_T e^{q_T} \Big|^r = 0.
\end{align}
Therefore applying the limit as $t \to T^-$ we 
get
\begin{align}
\notag
\lim_{t \to T^-}
\| u(x,t) - v_T (x) \|_{L^r}^r
&=
0.
\qed
\end{align}

Now that we have successfully established pointwise and $L^r$ convergence, we will use these results
to establish   
$H^s$ by using the Dominated Convergence Theorem.

\begin{theorem}
[Convergence in $H^s$  - Periodic version of Theorem \ref{nu-hs}]
\label{per-nu-hs} 
For $s < 5/4$, our antipeakon-peakon solution $u$ 
converges to $v_T$ in $H^s$, i.e.  
\begin{align}
\lim_{t \to T} \| u(x,t) - v_T(x) \|_{H^s} = 0.
\end{align}
\end{theorem}

{\bf Proof.}
From the definition of the $H^s$ norm of $u(x,t) - v_T(x)$,
we have
\begin{align}
\notag
\|u(x,t) - v_T(x) \|_{H^s}^2
&=
\sum_{n \in \zz}
(1 + n^2)^s
| \widehat{u}(n, t) - \widehat{v_T}(n)|^2.
\end{align}
Our objective, when taking the limit,  will be to move the limit inside of the integral.  Thus, 
we begin by examining the limit of the summand.
As we have already established
the convergence of $u$ to $v_T$  in $L^1$, via Proposition \ref{per-Lr-conv}, we see that
the inequality
\begin{align}
\notag
| \widehat{u}(n, t) - \widehat{u}(n, T)|
\le
\| u(x,t) - u(x,T) \|_{L^1},
\end{align}
implies that
\begin{align}
\notag
\lim_{t \to T^-}
(1 + n^2)^s
| \widehat{u}(n, t) - \widehat{v_T}(n)|^2 
&=
0.
\end{align}
Next, we will define the bounding
sequences that will allow us to apply the generalized Dominated Convergence Theorem.  We set
\begin{align}
\notag
f_t (n) \doteq
(1 + n^2)^s
|\widehat{u}(n, t) - \widehat{v_T}(n)|^2
&\le
4 (1+n^2)^s 
\Big(
|\widehat{u}(n, t)|^2 + |\widehat{v_T}(n)|^2  \Big) \doteq g_t (n). 
\end{align}
We need to establish that the $g_t$'s, have a pointwise limit $g$, i.e. for each $n \in \zz$, 
$\lim_{t \to T^-} g_t(n) = g(n)$.  The most
obvious candidate for $g$ is
\begin{align}
\notag
g(n)
=
8 (1 + n^2)^s |\widehat{v}_T(n)|^2
\quad \text{where} \quad
\sum_{n \in \zz} g(n)  =
32 \cosh^2 \pi \cdot c_s w_T.
\end{align}
Indeed, we have using the laws of limits, $\widehat{u}(n,t) \to  \widehat{v}_T(n)$
for each $n$ implies $g_t \to g$ for each $n$.
To finish satisfying the hypotheses of the generalized Dominated Convergence theorem, we must now
establish the sum properties of the $g_t$'s.
We have
\begin{align}
\notag
\lim_{t \to T}
\sum_{n \in \zz} g_t (n)  
&=
32 \sinh^2 \pi \cdot c_s w_T.
\end{align}
 where the left limit uses the $5/4$-hypothesis \eqref{limits}.
We now see that the hypotheses for the generalized Dominated Convergence Theorem are satisfied.
Thus, we can conclude that as $t \to T$, we have $u(x,t) \to v_T$ in $H^s$.
\qed

\vspace{.1in}

{\bf Proof of Theorem \ref{per-NE-non-uniqueness}.} 
Translating the NE  1-peakon solution 
by $x_0$ and choosing the minus sign we obtain the following antipeakon solution for NE
\begin{equation}
\notag
v(x,t)= - \sqrt{c}\,  \cosh ([x-x_0 -ct]_p - \pi), \quad \text{for any} \quad c>0
\quad
\text{and} \quad  x_0\in \tor.
\end{equation}
As in the real-line case, we choose
$c= w_T^2$ and $x_0 = q_T - w^2_T T,$
and obtain the antipeakon solution
$
v(x,t) = w_T  \cosh ([x  - x_0 - w_T^2 t] - \pi ).
$
Since at $t=T$ we have 
\begin{equation}
\notag
v(x,T)= w_T \cosh([x - q_T]_p - \pi)=u(x, T),
\end{equation}
we see that we have constructed two different NE solutions,
which belong in $H^s$, $s<5/4$, and agree at $t=T$.  
From here, a  change of 
variables can recast these two solutions as stemming from the same initial
data at time $t=0$.
This scenario proves failure of uniqueness in this range
of Sobolev spaces.  \, \qed

\vskip0.1in
{\bf The case $s=5/4$:} The argument for ill-posedness 
in this case is precisely the same 
as that in the non-periodic case.
This result completes the proof of both of Theorems
\ref{NE-2peak-slns} and \ref{NE-ILP}
in the periodic case.

%
%
%
%
%
%

\vskip0.1in

\noindent
{\bf Acknowledgements.} This work was partially supported by a grant from the Simons Foundation (\#246116 and \#524469 to Alex Himonas).
The first author thanks the Department
of   Mathematics of  the University of Chicago
for its hospitality during his 2015-2016 academic leave,
 where most of this  work was accomplished.
 The third author was supported in part by 
 NSF grants DMS-1265429 and DMS-1463746.
%

%
%
%
%
%
%
%
%
%

%

\vskip0.2in
\begin{minipage}[b]{9 cm}
   A. Alexandrou Himonas ({\it Corresponding author})  \\
    Department of Mathematics \\ 
     University of Notre Dame\\
 Notre Dame, IN 46556\\
      E-mail: {\it himonas.1$@$nd.edu}
\end{minipage}
\hfill
\begin{minipage}[b]{7 cm}
\noindent{Curtis Holliman}\\
Department of Mathematics\\
The Catholic University of America\\
Washington, DC 20064\\
E-mail: {\it holliman@cua.edu}
\end{minipage}
\vskip0.1in
  Carlos  Kenig\\
    Department of Mathematics \\  
 University of Chicago\\
5734 University Avenue\\
Chicago, IL 60637-1514\\
      E-mail: {\it cek@math.uchicago.edu}
%
%
 %
%
%
%
%
%
%
%
\end{document}